\documentclass[11pt,a4paper]{article}
\usepackage{geometry}
\usepackage{graphicx}
\usepackage{bm,array}
\usepackage{comment}
\usepackage{caption}
\usepackage{subcaption}
 \usepackage{makecell}
 \usepackage{multirow}
 \usepackage{float}

\usepackage[normalem]{ulem}

\usepackage[pdfpagemode={UseOutlines},bookmarks=true,bookmarksopen=true,
bookmarksopenlevel=0,bookmarksnumbered=true,hypertexnames=false,
colorlinks,linkcolor={blue},citecolor={blue},urlcolor={red},
pdfstartview={FitV},unicode,breaklinks=true]{hyperref}
\hypersetup{urlcolor=blue, colorlinks=true}

\usepackage{setspace}
\usepackage{vmargin}
\setmarginsrb{ 1.0in}  
{ 0.6in}  
{ 1.0in}  
{ 0.8in}  
{  20pt}  
{0.25in}  
{9pt}  
{ 0.3in}  

\usepackage{authblk}

\usepackage{amsmath}
\usepackage{amsfonts}
\usepackage{amssymb}
\usepackage[round, sort,comma,authoryear]{natbib}

\newtheorem{proposition}{Proposition}
\newtheorem{remark}{Remark}

\newenvironment{proof}[1][Proof]{\noindent\textbf{#1.} }{\ \rule{0.5em}{0.5em}}

\usepackage{multirow}
\usepackage{hhline}
\usepackage{footnote}
\usepackage[ruled,vlined]{algorithm2e}
\usepackage{algorithmic}

\newcommand{\transpose}{{\mbox{\tiny T}}}

\newcommand{\cV}{{\mathcal{V}}}
\newcommand{\cD}{{\mathcal{D}}}
\newcommand{\cS}{{\mathcal{S}}}

\newcommand{\cL}{{\mathcal{L}}}
\newcommand{\cX}{{\mathcal{X}}}

\newcommand{\cF}{{\mathcal{F}}}

\newcommand{\cY}{{\mathcal{Y}}}

\newcommand{\bC}{\textbf{C}}

\newcommand{\bx}{\textbf{x}}
\newcommand{\by}{\textbf{y}}
\newcommand{\bs}{\textbf{s}}
\newcommand{\bd}{\textbf{d}}

\newcommand{\bu}{\textbf{u}}
\newcommand{\bc}{\textbf{c}}

\newcommand{\bA}{\textbf{A}}

\newcommand{\bV}{\textbf{V}}

\newcommand{\ba}{\textbf{a}}
\newcommand{\bv}{\textbf{v}}
\newcommand{\bz}{\textbf{z}}
\newcommand{\bb}{\textbf{b}}

\newcommand{\bw}{\textbf{w}}
\newcommand{\bap}{\pmb{\alpha}} 
\newcommand{\bbt}{\pmb{\beta}}

\newcommand{\btheta}{\pmb{\theta}}
\newcommand{\bgamma}{\pmb{\gamma}}

\newcommand{\bbR}{\mathbb{R}}

\newcommand\myeq{\stackrel{\mathclap{\tiny\mbox{def}}}{=}}

\usepackage{mathtools}

\newif\ifnotes\notestrue
\def\mtien#1{{\color{magenta}{#1}}}

\def\htien#1{}

\newcommand{\ARUM}{\textsc{ARUM}}
\newcommand{\MRUM}{\textsc{MRUM}}

\usepackage{xcolor}

\begin{document}







\newcolumntype{C}{>{\centering\arraybackslash}p{4em}}

\title{\textbf{Joint  Location and Cost Planning in Maximum Capture Facility Location under Random Utilities}}
\author[1]{Ngan Ha Duong}
\author[1]{Tien Thanh Dam}
\author[1]{Thuy Anh Ta}
\author[2]{Tien Mai}
\affil[1]{\it\small 
ORLab,  Department of Computer Science, Phenikaa University}
\affil[2]{\it\small
School of Computing and Information Systems, Singapore Management University}

\maketitle


\begin{abstract}
We study a joint facility location and cost planning problem in a competitive market under random utility maximization (RUM) models. The objective is to locate new facilities and make decisions on the costs (or budgets) to spend on the new facilities, aiming to maximize an expected captured customer demand, assuming that  customers choose a  facility among all available facilities according to a RUM model. 
We examine two RUM frameworks in the discrete choice literature, namely, the additive and multiplicative RUM. While  the former has been widely used in  facility location problems, we are the first to explore the latter in the context. 
We numerically show that the two RUM frameworks can well approximate each other in the context of the cost optimization problem. In addition, we show that, under the additive RUM framework, the resultant cost optimization problem becomes highly non-convex  and may have several local optima. In contrast, the use of the multiplicative RUM brings several advantages to the competitive facility location problem. For instance,  
the cost optimization problem under the multiplicative RUM can be solved efficiently by a general convex optimization solver or can be reformulated as a conic quadratic program and handled by a conic solver available in some off-the-shelf solvers  such as CPLEX or GUROBI. Furthermore, we consider a joint location and cost optimization problem under the multiplicative RUM and propose three approaches to solve the problem, namely,  an equivalent conic reformulation, a multi-cut outer-approximation algorithm, and a local search heuristic. We provide numerical experiments based on synthetic instances of various sizes  to evaluate the performances of the proposed algorithms in solving  the cost optimization, and the joint location and cost optimization problems.

\end{abstract}

{\bf Keywords:}  
Competitive facility location, maximum capture, joint location and cost optimization, 
multiplicative random utility maximization, convex optimization, conic programming, outer-approximation, local search heuristic 


%


\section{Introduction}

Facility location has been one of the most  fundamental problems in decision-making for modern transportation and logistics systems. In the most classical setting, this is the problem of  determining a subset of locations from a given set of location candidates to set up new facilities, aiming at maximizing a profit function (e.g. expected customer demand, or expected revenue) or minimizing a cost function (e.g., operational or transportation cost). Customer demand plays a critical role in facility location problems as it is a driving factor in determining where to open facilities. In this work, we focus on a class of competitive facility location problems where customers' demand is modeled and predicted  by a random utility maximization (RUM) model \citep{BenaHans02,Ljubic2018outer,mai2020multicut,Dam2021submodularity}. In this class of problems, it is assumed that customers choose among different facilities by maximizing  utilities that they assign to the available facilities. Such utilities are typically
functions of facility attributes/features, e.g., quality of services or infrastructure of the facilities, or transportation costs. The use of the RUM framework  in the context is well supported by the popularity and advantages  of RUM models in modeling and predicting people's choice behavior \citep{BenALerm85}. To the best of our knowledge, existing RUM-based facility location studies only focus on making location decisions, ignoring the fact that customers' demand would be significantly affected by other factors such as the opening/operating costs spent on the new facilities \citep{aboolian2007competitive,ciari2014modeling}.  We address this limitation in this paper.  

It is well-known that facility location and cost decisions would significantly impact customer demand, thus determining the success of the business, especially in a competitive market where other players are competing for the market share. For instance, in the context of car-sharing businesses where customers have the options of using their own vehicles or renting cars from different car-sharing services,  \cite{jorge2013carsharing,ciari2014modeling} show  that customer demand  is mainly affected by the distance to the car-rental stations and the service availability at these stations, and \cite{shaheen1998carsharing} indicate that when customers observe more vehicle availability and usage conveniences in the newly opened rental stations, their demand for the service  would increase. Intuitively, customer demand would be critically affected by both the service availability within
their neighborhood (i.e. locations  of the service) and the quality of the service, which is generally dependent on the costs spent on the facilities. This is our main motivation to jointly consider  the  location and cost planning in the context of competitive facility location.  

In our joint facility and cost optimization problem, the  objective is to maximize an expected captured customer demand by making decisions on choosing locations to open new facilities and distributing the available budget to the  new facilities. Since the objective is to maximize a captured demand function, such a  problem is also referred to as the maximum capture problem (MCP).
To model the impact of the costs on customer behavior, we include the costs into the customers' utility functions, taking into consideration the fact that spending a higher budget on a facility would make this facility more attractive to customers and increase the probability that the facility will be chosen. 
We then formulate the problem as an optimization problem with binary (for the choices of locations) and continuous (for the costs) variables.  Due to the nonlinear structures of the RUM models,  the resulting location and cost optimization problem becomes highly nonlinear  and challenging to handle.  We then focus on the question of how to efficiently solve the proposed problems under large-scale settings.
We will state our main contributions in the following, but before that, we note that, throughout the paper, when we say ``choice-based decision-making problem'', we are  referring to a decision-making problem where customer/user behavior is predicted by a discrete choice model, e.g., facility location or assortment optimization under a RUM discrete choice model. 

\textbf{Contributions:}
    We consider both the additive  RUM (ARUM) \citep{McFa78} and multiplicative RUM (MRUM) \citep{FosgBier09} frameworks, where  the former assumes that a random utility is a  sum of a deterministic and a random error term and  the latter framework relies on the setting that the random utilities have a multiplicative structure, i.e., a random utility is a product of a deterministic and a random error term. 
    While the former has been  popularly used in facility location (and other choice-based decision-making problems), we are the first to consider the latter in the context. 
    We first provide a numerical study on the relation between the ARUM and MRUM frameworks in the context of the MCP. More precisely, we experimentally show that the choice probability functions and the MCP objective functions under the ARUM and MRUM can be well approximated by each other, supporting the use of the MRUM framework even when ARUM is the ``ground-truth'' one.
    
    We further show that while the cost optimization problem under ARUM is highly non-convex and may have several local optimal solutions, the cost optimization problem under the MRUM framework can be solved efficiently by convex optimization. In addition, we show that the cost optimization problem under MRUM can be reformulated as a conic quadratic program, which can be efficiently handled by an off-the-shelff solver (e.g. CPLEX or GUROBI). These findings indicate the fact that the MRUM framework would bring great advantages in the context of choice-based facility location, as compared to the classical ARUM framework. 
    
    We then consider the joint location and cost optimization problem under MRUM and propose three solution approaches to solve the problem. First, we propose an equivalent conic reformulation that allows us to exactly solve the problem by an off-the shelf solver such as CPLEX or GUROBI \citep[see][for instance]{bonami2015recent}. Second, we show that the joint problem can be also handled exactly by a multicut outer-approximation algorithm \citep{Bonami2011_BB_MIP,hijazi2014outer,mai2020multicut}. 
    Third, motivated by the fact that some recent studies have shown that a local search heuristic would achieve state-of-the-art performance in the context \citep{Dam2021submodularity}, we develop such a local search procedure to solve our joint problem. We summarize our findings and solution methods in Table \ref{tab:summary} below. 
    
\begin{table}[htb]
\centering
\begin{tabular}{l|l|l}
Problem                                                                        & MCP under ARUM                                                                   & \textbf{MCP under MRUM}                                                                                                                                                                              \\ 
\hline
Cost optimization                                                              & \begin{tabular}[c]{@{}l@{}}Highly non-convex,\\intractable to solve\end{tabular} & \begin{tabular}[c]{@{}l@{}}\textbf{Exact methods:}\\~ ~- Convex optimization solver\\~ ~- Conic reformulation~\end{tabular}                                                                            \\ 
\hline
\begin{tabular}[c]{@{}l@{}}Joint location \\and cost optimization\end{tabular} & Intractable                                                                      & \begin{tabular}[c]{@{}l@{}}\textbf{\textbf{Exact methods:}}\\~ ~- Conic reformulation\\~ ~- Outer Approximation\\\textbf{\textbf{\textbf{\textbf{Heuristic method:}}}}\\~ ~- Local search\end{tabular}  \\
\hline
\end{tabular}
\caption{Solution methods.}
\label{tab:summary}
\end{table}
   
    We provide numerical experiments based on synthetic instances of various sizes to validate the performances of our approaches, for both the cost optimization and the joint location and cost optimization problems under the MRUM framework. Our numerical results clearly demonstrate the superiority of the conic reformulation approach for the cost optimization problem, and the multicut outer-approximation algorithm for the joint location and cost optimization problem, especially for large-scale problem instances.

In this paper, we focus on multinomial logit (MNL) based models, i.e., the random terms are independent and identically distributed (i.i.d.) \citep{Trai03}. This is the most popular model in the discrete choice modeling and choice-based decision-making literature. Other choice models taking into consideration the correlation between choice alternatives exist, e.g., the nested or cross-nested logit models \citep{McFaTrai00,Trai03}, but
 the MCP under such models are generally intractable to solve \citep{mai2020multicut,Dam2021submodularity} and 
 we  keep them for future work. 

\textbf{Literature review:}
Our work contributes to  the line of research on the maximum capture facility location problem under random utility maximization discrete choice models. In this context, \cite{BenaHans02} seems to be the first to formulate the MCP under the ARUM-based multinomial logit model (MNL) and  propose methods based on mixed-integer linear programming (MILP) and variable neighborhood search. Afterward, some alternative MILP reformulations have been proposed by \cite{FLO_Zhang2012impact} and \cite{FLO_Hasse2009MIP}. \cite{FLO_Haase2014comparison}  provide a comparison of existing MILP reformulation models and conclude that the MILP model from \cite{FLO_Hasse2009MIP} gives the best performance. \cite{FLO_Freire2016branch} strengthen the MILP reformulation
of \cite{FLO_Hasse2009MIP} using a branch-and-bound algorithm with some tight inequalities.  \cite{Lin2021branch,Lin2021generalized} further improve the MILP approach using Benders decomposition. \cite{Ljubic2018outer} propose a branch-and-cut method combining 
outer-approximation and submodular cuts, and \cite{mai2020multicut} propose a multicut outer-approximation algorithm to efficiently solve large-scale instances. Recently, \citep{Dam2021submodularity} study the MCP under the generalized extreme value (GEV) family of models and propose a local search algorithm with guarantees. To the best of our knowledge, existing works in the context of the MCP under RUM  only focus on the selection of locations, yielding optimization problems that only involve binary variables. In contrast, our joint location and cost optimization problem involves both binary and continuous variables, thus existing \textit{exact} combinatorial optimization algorithms do not directly apply.  

It is worth mentioning that there have been papers proposing to use simulation to approximate the  choice model, yielding mixed-integer linear programs (MILP) for the approximated choice-based optimization problems  \citep{paneque2021integrating,lamontagne2022optimising}.
This approach is general as any RUM based choice model can be approximated by simulation, but is limited in the sense that  a large number of samples would be required to deliver a good approximation of the choice model and  the size of the MILP grows quickly as the number of samples increases. In contrast, our approach is exact and the sizes of our MILP and conic reformulations are remarkably smaller as they are only proportional to the number of locations.  
Our work also belongs to the broad literature of facility location  \citep[see][for a review]{Laporte2015introduction}, noting that decisions on both locations and costs have been widely considered in the broad context
\citep{aboolian2007competitive,berman2009modeling}.

The literature on random utility maximization models covers both the ARUM and MRUM frameworks. While the ARUM framework has been widely studied since the 1908s and already has as many successful applications \citep{McFa78,McFa81,Trai03}. The MRUM framework appears in the literature since the 2000s and seems to be less popular. \cite{MRUM_brilon2002multiplicative} seem to be the first to mention the multiplicative structure and since then this framework has been studied and applied by various researchers \citep{FosgBier09, MRUM_hess2018revisiting,MRUM_borjesson2012valuations,MRUM_borjesson2014experiences,MRUM_de2014new,MRUM_fosgerau2006investigating}. In particular, it has been shown that  changing from the additive to the multiplicative structure (i.e, from ARUM to MRUM) would  lead to a large improvement in model fit, sometimes larger than that gained from introducing
random coefficients in the deterministic terms of the utilities (i.e., the mixed MNL model  \citep{McFaTrai00}) \citep[see][for instance]{FosgBier09}. As far as we know, in the context of the choice-based facility location and other choice-based decision-making  problems such as assortment and price optimization \citep{Talluri2004revenue,rusmevichientong2010dynamic,rusmevichientong2014assortment}, our work  is the first time the MRUM framework is explored. 

Since we study a joint location and cost optimization problem, our work closely relates to joint assortment and pricing optimization problems \citep{wang2012capacitated,chen2020capacitated,miao2021dynamic}, where the relation
to the MCP problem stems from the analogy of set of products to
set of locations and of continuous prices to continuous costs to spend on the new facilities. In particular, the objective function in the MCP is a sum of ratios, thus our problem  share a similar structure with a joint assortment and price optimization problem under the mixed MNL model \citep{McFaTrai00}, for which \cite{li2019product} show that the price optimization itself is highly non-convex and challenging to solve. Later on in the paper, we will also show that the  cost optimization MCP under ARUM is also highly non-convex and  may have several local optima. Again, we highlight the fact that, as far as we know, existing assortment and/or price optimization work only employs the ARUM framework. Thus, our work here would encourage the use of the MRUM counterpart in assortment and/or price optimization problems.

\noindent
\textbf{Paper outline:} Our paper is organized as follows. Section \ref{sec:prob-formulation} presents the problem formulations where we introduce the formulations of the joint location and cost optimization problems under both ARUM and MRUM frameworks. Section \ref{sec:relation RUM} analyse the relation between the ARUM and MRUM frameworks in the context of the MCP. Section \ref{sec:methods} presents our solution methods to solve the cost optimization and joint location and cost optimization problems under MRUM. Section \ref{sec:experiments} shows our numerical evaluations, and Section \ref{sec:concl} concludes.  

\noindent
\textbf{Notation:}
Boldface characters represent matrices (or vectors), and $a_i$ denotes the $i$-th element of vector $\ba$. We use $[m]$, for any $m\in \mathbb{N}$, to denote the set $\{1,\ldots,m\}$. 


\section{Problem Formulations}
\label{sec:prob-formulation}
In this section, we present the additive  and multiplicative RUM models and introduce the formulations of the joint facility location and cost optimization problems under these frameworks.
\subsection{Additive and Multiplicative RUM models}
The discrete choice framework \citep{Trai03} assumes that each individual (decision-maker) $n$ associates an utility $U_{ni}$ with each alternative $i$ in a choice set $S_n$. 
It is then assumed that each utility $U_{ni}$ consists of two parts: a deterministic part $V_{ni}$ that contains some observed attributes of the corresponding alternative and the individual $n$, and a random term $\epsilon_{ni}$ that is unknown to the analyst. Different assumptions
 can be made on the random terms, leading  to different types of discrete choice models. In general, a
linear-in-parameters formula is used for the deterministic part of the utility, i.e, $V_{ni} = \bbt^{\transpose} \ba_{ni}$, where $\bbt$ is a vector of parameters that  can be estimated using historical observations and $\ba_{ni}$ is a vector of attributes
of alternative $i$ as observed by individual $n$. Under the RUM framework, the probability that individual $n$ selects item $i\in S_n$ can be computed as  $P(U_{ni}\geq U_{nj},\;\forall j\in S_n)$, i.e., the individual chooses an item of the highest utility. 

Under the additive RUM (ARUM) framework \citep{McFa78,Trai03}, the deterministic parts and random terms have an additive form, $U_{in} = V_{
ni}+\epsilon_{ni}$. Under the popular multinomial logit (MNL) model where  $\epsilon_{ni}$ are i.i.d  and  follow the standard Gumbel distribution, the choice probabilities have the following form 
\[
P(i|S_n) = \frac{e^{V_{ni}}}{\sum_{j\in S_n} e^{V_{nj}}},\;\forall i\in S_n.
\]
The MNL is widely used in both demand modeling  and decision-making \citep{Trai03,Talluri2004revenue} due to its simple structure. The model however retains  the independence  from irrelevant alternatives (IIA) property which may be not holds in many situations. Efforts have been made to overcome this shortcoming, resulting in several discrete choice models, i.e., the nested logit  \citep{BenA73}, cross nested logit \citep{VovsBekh98}, network multivariate extreme value (MEV)   \citep{DalyBier06}, or mixed logit \citep{McFa81} models.  

The multiplicative RUM framework \citep{MRUM_brilon2002multiplicative,FosgBier09,MRUM_borjesson2012valuations} replaces the additive assumption in the  ARUM framework by its multiplicative couter part, i.e., $U_{ni} = V_{ni}\epsilon_{ni}$, where $V_{ni}$ and the random terms $\epsilon_{ni}$ are additionally assumed to be positive to make $\ln V_{ni}$ and $\ln\epsilon_{ni}$ valid to derive the choice probabilities. The choice probability that individual $n$ selects item $i$ becomes $P(V_{ni}\epsilon_{ni}\geq V_{nj}\epsilon_{nj}\;\forall j\in S_n)$. Since  the logarithm is a strictly increasing function, the choice probabilities can be further written as 
\[
P(V_{ni}\epsilon_{ni}\geq V_{nj}\epsilon_{nj}\;\forall j\in S_n) = P( \ln V_{ni} +\ln \epsilon_{ni} \geq \ln V_{nj} +\ln \epsilon_{nj},\;\forall j\in S_n).
\]
We then can impose  some similar assumptions on the random terms $\ln\epsilon_{nj}$, as in the ARUM framework, to get closed form choice probabilities. For instance, if $\ln \epsilon_{ni}$ are i.i.d  and follow the standard Gumbel distribution, the choice probabilities become 
\[
P(i|S_n) = \frac{e^{\ln V_{ni}}}{\sum_{j\in S_n}e^{\ln V_{nj}}} = \frac{V_{ni}}{\sum_{j \in S_n} V_{nj}}. 
\]
Here, we note that it is only necessary to assume that the signs of $V_{ni}$ are the same over $i\in S_n$. In the case that $V_{ni}$ are negative, we can write the choice probabilities  as $P((-V_{ni})(-\epsilon_{ni})\geq (-V_{nj})(-\epsilon_{nj}),\;\forall j\in S_n) = P(\ln (-V_{ni}) +\ln(-\epsilon_{ni})\geq \ln(-V_{nj}) + \ln(-\epsilon_{nj}),\;\forall j\in S_n)$, with the assumption that $-\epsilon_{ni}$ are positive. One then can assume that $\ln(-\epsilon_{ni})$ are i.i.d Gumbel distributed to obtained a MNL-based model $P(i|S_n)  = (-V_{ni})/\left(\sum_{j\in S_n}(-V_{nj})\right)$.
Similarly  to the ARUM framework, other assumptions can be made for the random terms ($\ln \epsilon_{ni}$ or $\ln (-\epsilon_{ni})$) to obtain models with nested, cross nested, or mixed structures. 

\subsection{Joint Location and Cost Optimization MCP}
In the context of the competitive maximum capture problem, a ``newcomer'' firm aims to locate new facilities in a competitive market where there are existing facilities from competitors. The firm's objective is to maximize the expected captured market share achieved by attracting customers to new facilities. 
To capture the customers’ demand, we suppose that a customer selects a facility according to a RUM discrete choice model (ARUM or MRUM). 
In this context, each customer $n$ associates each new facility $i$ with a random utility $U_{ni}$  and choices are made by maximizing his/her utilities. 
The firm's decisions are a set of locations to setup new facilities and costs to spend on new facilities. In the following, we describe in detail our problem formulations under both ARUM and MRUM MNL-based models. 

Assume that there are $m$ possible locations to locate new facilities. Let  us denote by $\cV= [m]$ the set of all possible locations and by $I$ the set of all  customer zones. Depending on the application, $I$ can be a set of customer types categorized by customers' characteristics (e.g, age or income), or can be a set of geographical zones where customers are located. Assume that there are $q_n$  customers for each zone $n \in I$. Let $V_{ni}$ be the deterministic utility that customer zone $n\in I$ assign to  a facility at location $i\in [m]$ when making choice decisions, and $\bV^n$ is a vector of size $m$ with entries  $V_{ni}$, for all $n\in I$.  
Given any subset $S$ of possible locations i.e. $S\subset [m]$ under the additive MNL model, the choice probabilities of customer zone $n\in I$ are given as
\[
P^{\ARUM}(i|S,\bV^n) = \frac{e^{V_{ni}}}{U^c_n+ \sum_{i'\in S} e^{V_{ni'}}},  
\]
where $U^c_n$ represents the customer zone n’s utility of choosing a competitor’s facility. 
When the facility  planing and cost optimization problems are considered simultaneously, we write the utilities $V_{ni}$ as  $V_{ni} = a_{ni} x_i + b_{ni}$, where $a_{ni}$ is a parameter representing the sensitivity of customer zone $n$ w.r.t. the cost spent on facility $i$, $x_i$ is the cost to spend on facility $i$, and $b_{ni}$ is another factor that affects customer choice decisions, e.g., distance or parking availability.
The parameters $a_{ni}$ play an important role in our model, in the sense if a customer zone $n$ is located far away from location $i$, the investing cost of facility $i$ would not make a significant impact on the choices of the customers of this zone, and on the other hand if location $i$ is close to the location of customer zone $n$, then adding more budget to the facility setup would make a huge impact on the attractiveness of the facility. In fact, it is to be expected that larger $x_i$ (more cost to spend on facility $i$) will make the facility more attractive to customers, thus $a_{ni}$ should take positive values. 
The joint facility location and cost optimization problem can be formulated as 
\begin{align}
     \max_{S\in\cS,~\bx} &\left\{f^{\ARUM}(S,\bx) =\sum_{n\in I} q_n\sum_{i\in S} \frac{e^{a_{ni} x_i +b_{ni} }}{U^c_n+ \sum_{i'\in S} e^{a_{ni'}x_{i'} + b_{ni'}}} \right\} \nonumber \\
    \mbox{subject to} &\quad \sum_{i \in S} x_i \leq C, \nonumber\\
    &\quad x_i \in [L_i,U_i], \; \forall i \in S \label{eq:ctr-bounds}\\
    &\quad \bx \in \bbR^{|S|}_+, \nonumber
\end{align}

where $C$ is the maximum budget to spend on the new facilities, $L_i$  and $U_i$ 
are the lower and upper bounds for the cost on  facility $i\in [m]$, $\cS$ is the set of all subsets of [m]. It is required that $\sum_{i\in S}L_i \leq C$ for some $S\in\cS$ to ensure feasibility. More general linear constraints on the costs $\bx$ can be included, but we only consider a budget constraint $\sum_{i\in S}x_i\leq C$ and the bound constraints for the sake of simplicity. We can represent set $S$ by binary variables vector with $y_i= 1$ if set S consists of location $i$ and $y_i = 0$ otherwise. Then, we formulate the MCP under ARUM as a binary-continuous nonlinear program
\begin{align}
     \max_{\by,\bx} &\left\{f^{\ARUM}(\by,\bx) =\sum_{n\in I} q_n \frac{ \sum_{i\in[m]} y_i e^{a_{ni} x_i + b_{ni} }}{U^c_n+ \sum_{i\in [m]} y_ie^{a_{ni}x_{i} + b_{ni}}} \right\}   \label{prob:MCP-addtive}\tag{\sf MCP-ARUM}      \\
    \mbox{subject to} &\quad \sum_{i \in [m]} x_iy_i \leq C \label{eq:ctr-capacity}\\
    & \quad \sum_{i\in [m]} y_i \leq K \label{eq:ctr-capacity-M}\\
    &\quad x_i \leq y_i U_i,\;\forall i\in [m]\label{eq:ctr-x-y-j} \\
    &\quad x_i \geq  y_i L_i,\;\forall i\in [m]\label{eq:ctr-x-y-j-lb} \\
    &\quad \bx \in \bbR^{m}_+,\; \by\in \{0,1\}^m   \nonumber
\end{align}
where $K$ %
is the maximum number of facilities to be opened,
constraints \eqref{eq:ctr-x-y-j} is to ensure that if there is no facility at  location $i$, the corresponding cost should be zero and \eqref{eq:ctr-x-y-j-lb} is to force the cost $x_i$ to be at least $L_i$ if location $i$ is selected (i.e., $y_i = 1$). 

If the choice model is a multiplicative MNL model, the choice probabilities become
\[
P^{\MRUM}(i|S,\bV^n) = \frac{{V_{ni}}}{1 +  \sum_{i'\in S} {V_{ni'}}} = \frac{{a_{ni}x_i+b_{ni}}}{U^c_n +  \sum_{i'\in S} {a_{ni'}x_{i'}+b_{ni'}}},
\]
and the joint facility and cost optimization MCP can be formulated as
\begin{align}
     \max_{\by,\bx} &\left\{f^{\MRUM}(\by,\bx) =\sum_{n\in I} q_n \frac{\sum_{i\in [m]} y_i({a_{ni} x_i +b_{ni} })}{U^c_n+ \sum_{i\in [m]} y_i({a_{ni}x_{i} + b_{ni}})} \right\}  \label{prob:MCP-mulpli}\tag{\sf MCP-MRUM}  \\
    \mbox{subject to} &\quad \text{Constraints \eqref{eq:ctr-capacity}-\eqref{eq:ctr-capacity-M}-\eqref{eq:ctr-x-y-j}-\eqref{eq:ctr-x-y-j-lb}}, \nonumber\\
    &\quad \bx \in \bbR^{m}_+,\; \by\in \{0,1\}^m.  \nonumber
\end{align}

If the costs $\bx$ are fixed  and only the facility location problem is considered, then \eqref{prob:MCP-addtive} and \eqref{prob:MCP-mulpli} share the same structure, thus existing methods developed for the ARUM-based MCP problem \citep{FLO_Benati2002maximum,FLO_Haase2014comparison,mai2020multicut} can be applied to the MRUM-based problem. When the cost optimization is taken into consideration, the joint facility location and cost optimization problem under MRUM seems  to have a simpler structure, as exponential functions  are now replaced by linear ones. We discuss this in more detail in the next sections. 



\section{Relation between MRUM and ARUM Models }\label{sec:relation RUM}
We numerically analyse the relation between the ARUM and MRUM frameworks in the context of the MCP. We employ MNL-based models (i.e., the random terms are i.i.d extreme value Type I) for our analyses, setting our focus on the impact of the structures of the deterministic terms under the MRUM and ARUM on the choice probabilities and on the objective function of the MCP. Since the ARUM framework is well understood and widely studied, we aim to explore the question of whether the MRUM models and the MCP under MRUM are able to well approximate their counterparts under the ARUM.

\begin{figure}[htb] 
\centering
    \includegraphics[width=1.0\linewidth]{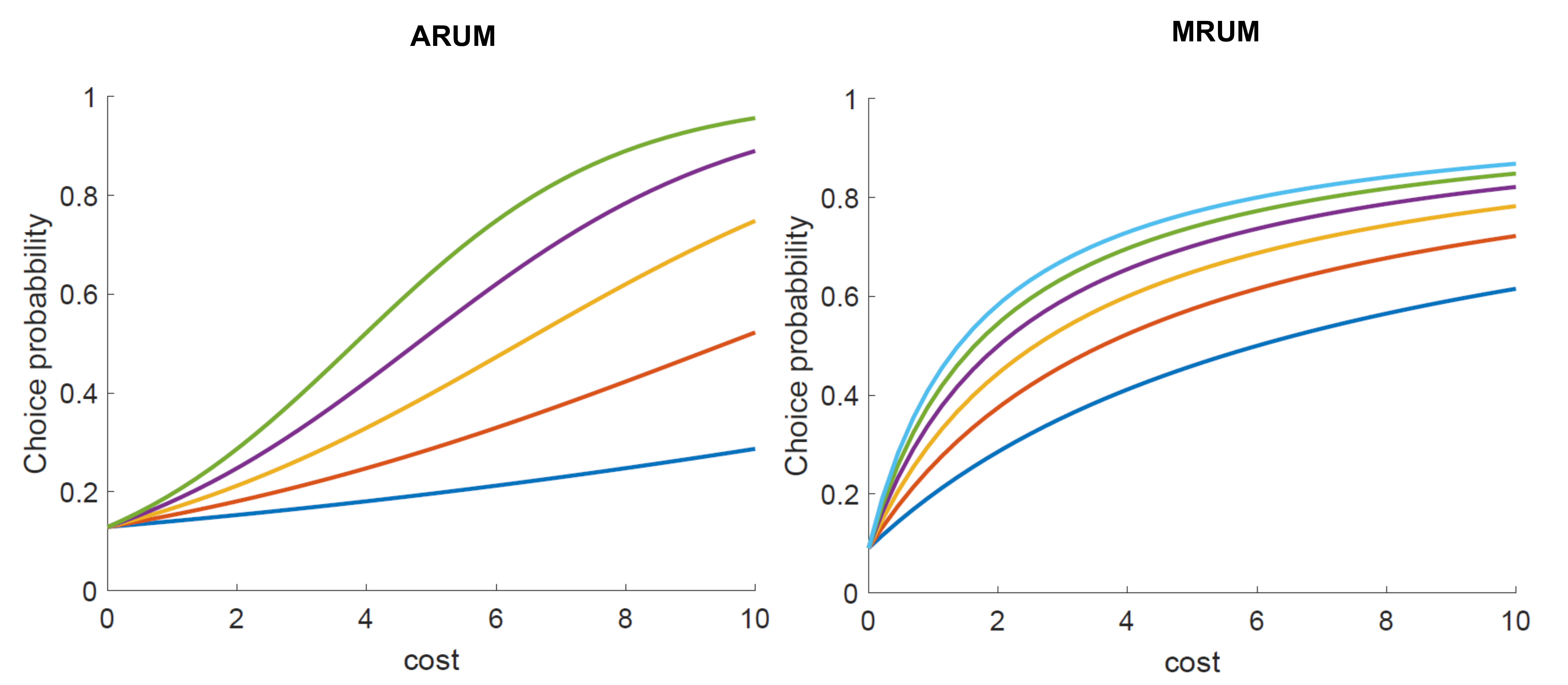} 
    \caption{Choice probability functions under ARUM and MRUM.} 
    \label{fig:prob-MCP-ARUM} 
\end{figure}

We first note that the choice probability of an alternative $j$ under ARUM is of the form $(\exp(a_jx_j+b_i))/(\exp(a_jx_j+b_j) + c)$, where  $a_j,b_j$ are some choice parameters, $x_j$ is the investment cost of $j$, and $c$ represents the utilities of the other alternatives. On the other hand, under MRUM,  the choice probability of an alternative $j$ under ARUM is of the form $((a'_jx_j+b'_i))/((a'_jx_j+b'_j) + c')$, where $a'_j,b'_j$ are choice parameters and  $c'$ is the sum of the utilities of the other  alternatives. We first explore the shapes of the choice probabilities, under the ARUM and MRUM, as functions of the cost. To this end, we select some sets of parameters $(a_j,b_j,c,a'_j,b'_j,c')$ and vary the cost $x_j$. The choice probabilities under the ARUM and MRUM are plotted in Figure \ref{fig:prob-MCP-ARUM}. It can be seen that the choice probability function under MRUM  has a concave shape, while it is not the case under the ARUM. This can be validated by seeing that function $((a'_jx_j+b'_i))/((a'_jx_j+b'_j) + c')$ is concave in $x_j$.
Moreover, under MRUM, the choice probability seems to increase faster when $x_j$ is small, but slower when  $x_j$ becomes large, as compared to the choice probability function under ARUM.

\begin{figure}[htb] 
\centering
    \includegraphics[width=1.0\linewidth]{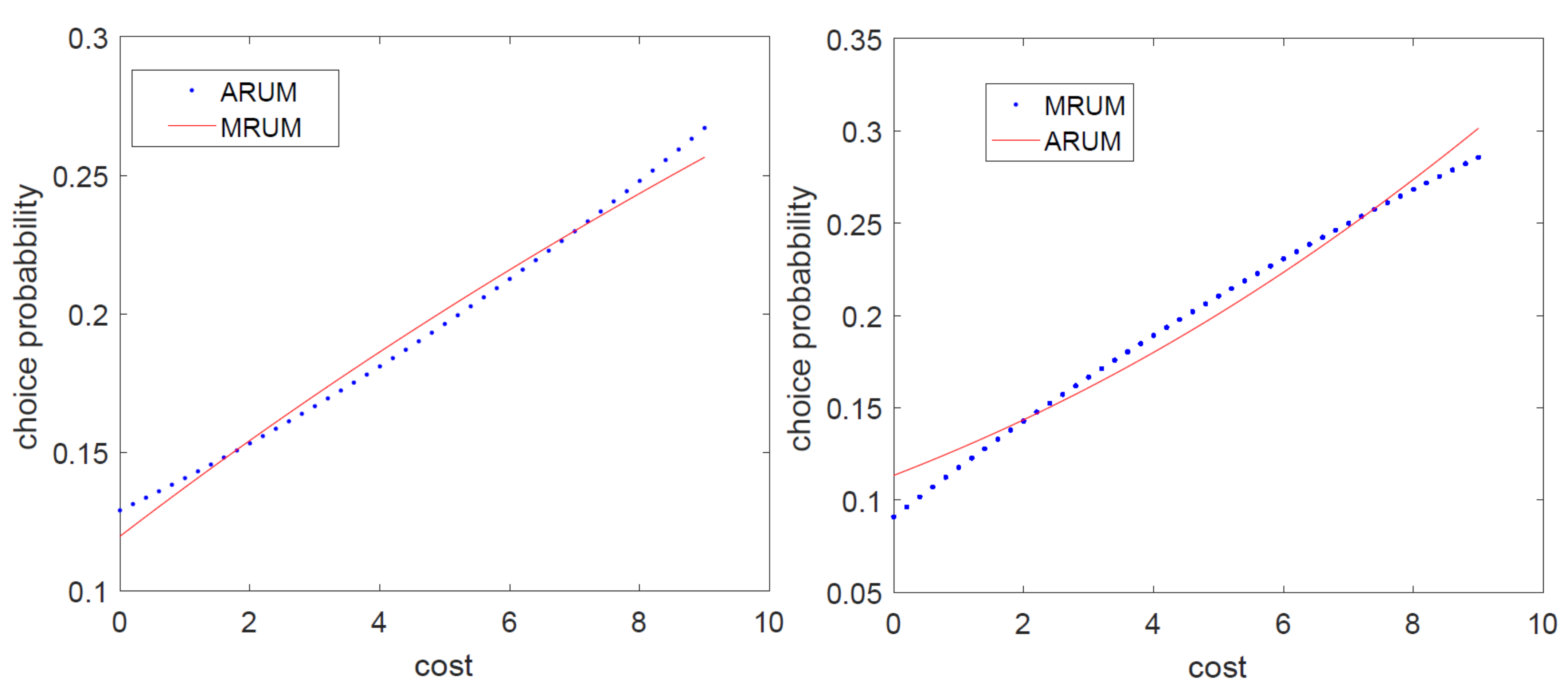} 
    \caption{Fitting choice probability functions yielded by the two frameworks in two-dimensional space.} 
    \label{fig:prob-MCP-ARUM} 
\end{figure} 
In the second experiment, we aim to explore whether the MRUM can well approximate the choice probabilities given by the ARUM, and vice versa. To this end, we choose a ARUM model with choice probabilities of the form $(\exp(a_jx_j+b_i))/(\exp(a_jx_j+b_j) + c)$. We then vary $x_j$ to get pairs of costs and choice probabilities $\{(x_j^1,p_j^1), (x_j^2,p_j^2)\ldots\}$. Next, we fit these data points with a MRUM choice probability function of the form $((a'_jx_j+b'_i))/((a'_jx_j+b'_j) + c')$, using Least-Squares regression. We do the same to fit an ARUM choice probability function with data points generated by a MRUM model. In the left figure of Figure \ref{fig:prob-MCP-ARUM} we generate data points using function $(\exp(0.1x_j-0.3))/(\exp(0.1x_j-0.3)+5)$  and get $((0.14x_j+0.82))/((0.14x_j+0.82) +6.02)$ after fitting, with a Root Mean Square Error (RMSE) of 0.0049 (about 2\%). On the right hand side of Figure \ref{fig:prob-MCP-ARUM} we generate the data points using $((x_j+3))/(((x_j+3))+30)$  and get $(\exp(0.14x_j+0.1))/(\exp(0.14x_j+0.1) +8.67)$ after fitting, yielding a RMSE of 0.0092 (about 4\%). In general, our experiment shows that, even-though the choice probability functions under ARUM and MRUM have different shapes, they can well approximate each other.  
\begin{figure}[htb] 
\centering
    \includegraphics[width=1.0\linewidth]{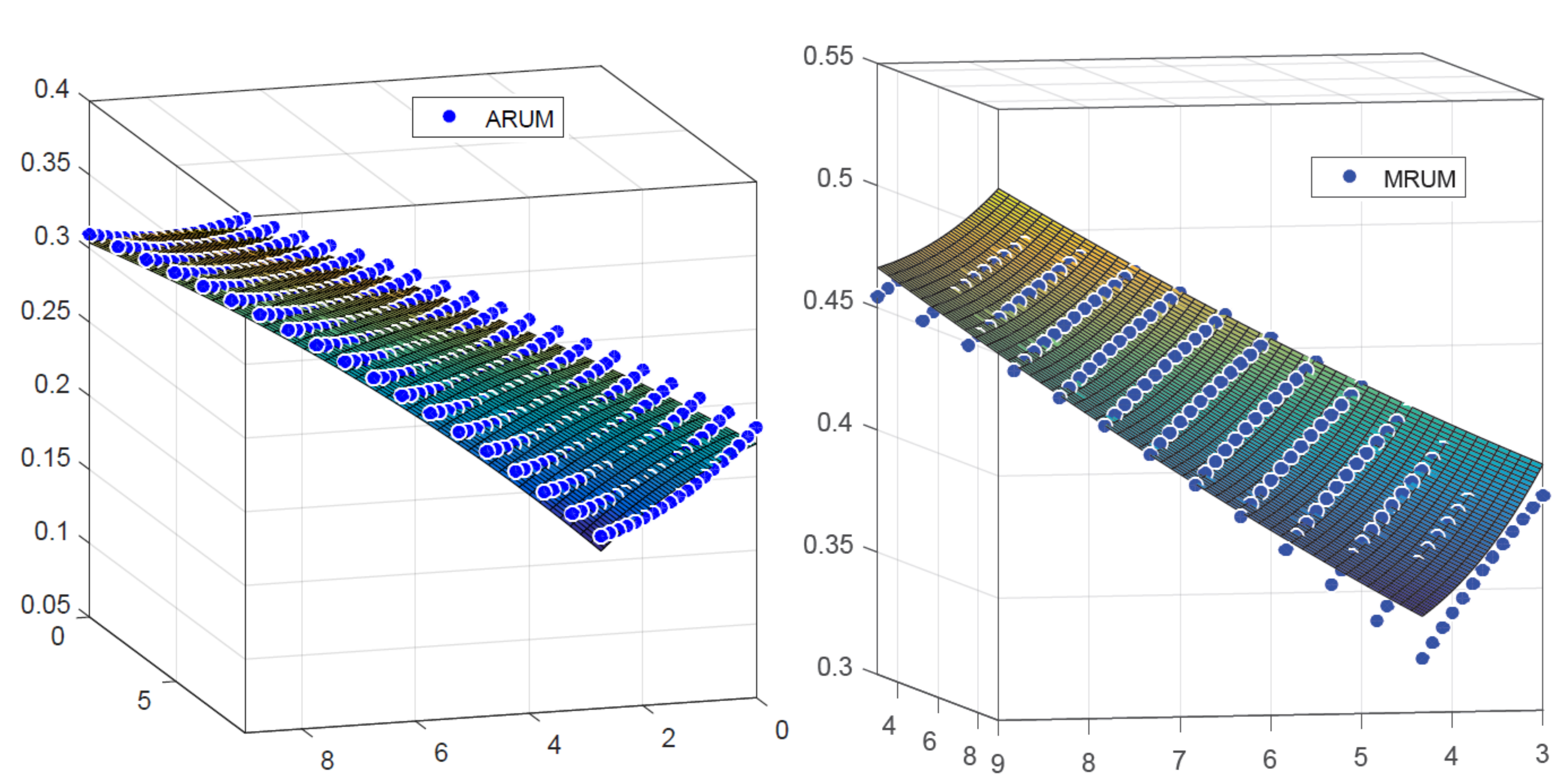} 
    \caption{Fitting choice probability functions yielded by the two frameworks in 3D space.} 
    \label{fig:prob-MCP-ARUM-3D} 
\end{figure}

We extend the second experiment to a three-dimensional (3D) space  to see how the MCP objective functions can approximate each other in higher dimensions. To this end, we select MCP objective functions  depending on two cost variable $x_i,x_j$ to generate data points. We plot the data points (blue dots)  and the fitted functions (surfaces) in Figure \ref{fig:prob-MCP-ARUM-3D}, where the left hand side figure show data points from the function $(\exp(0.1x_i-0.3)+(0.2x_j-0.3))/(\exp(0.1x_i-0.3)+(0.2x_j-0.3)+5)$ and a surface of the fitted function $((0.15x_i-0.24)+(0.25x_j-0.24))/((0.15x_i-0.24)+(0.25x_j-0.24)+6.25)$, with a RMSE of 0.0044. On the right hand side of Figure \ref{fig:prob-MCP-ARUM-3D} we plot data points from the function $((x_i+3)+(2x_j+1))/(((x_i+3)+(2x_j+1))+30)$ and a surface of the function $(\exp(0.16x_i-0.36)+\exp(0.13x_j+0.71))/(\exp(0.16x_i-0.36)+\exp(0.13x_j+0.71)+8.8)$, with a RMSE of 0.0054, which is just about 1\% of the average of the MCP objective values.
In general, the experiment shows that the MCP objective functions under the ARUM and MRUM can well approximate each other in a 3D space.

\begin{figure}[H] 
\centering
    \includegraphics[width=0.7\linewidth]{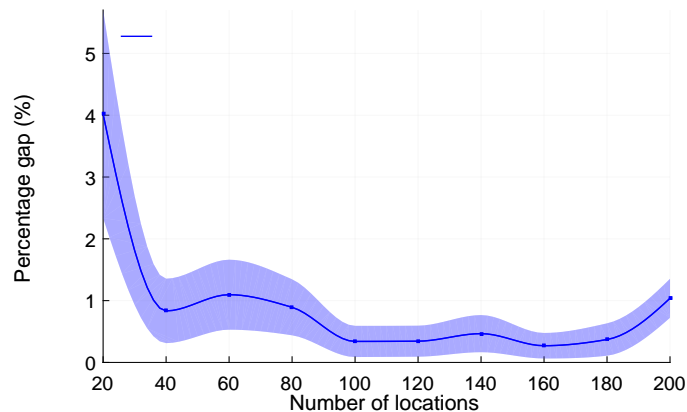} 
    \caption{Percentage gaps between the objective values of the MCP under the MRUM and fitted ARUM models. } 
    \label{fig:MCP-MRUM-ARUM} 
\end{figure}

In the fourth experiment, we further explore whether the MRUM framework can approximate the ARUM in the context of the MCP. We randomly create a ``ground-truth'' ARUM model and use it to generate 1000 choice observations. We then specify a MRUM model of the same utility specification and estimate it's parameters using the generated choice observations by maximum likelihood estimation. To compare the ARUM model and the estimated MRUM one, we randomly and uniformly generate 1000 samples of the  cost vector $\bx$ and compute the percentage gaps between the MCP objective functions under the two choice models. We vary the number of alternatives from 20 to 200 and plot the means (solid curve) and standard errors (shaded area) of the percentage gaps in Figure \ref{fig:MCP-MRUM-ARUM}. 
It can be seen that the percentage gap, on average, is only 4\% when $m=20$ and just about $1\%$ for the other cases, showing that even when the ``ground-truth'' choice model is an ARUM, one still can well approximate it with small gaps by a MRUM model, noting that the MCP under MRUM is more tractable to solve (see Table \ref{tab:summary} above).

\section{Solution Methods} \label{sec:methods}
In this section, we develop solution methods to solve the joint location and cost MCP. As will be shown later, the MCP under ARUM is highly non-convex and intractable to handle, we will set our focus on the  MCP under MRUM.
\subsection{Cost Optimization  MCP}
\label{sec:CP}
We first consider the cost optimization problem, 
assuming that the location variables $\by$ are fixed. We show that the cost  optimization problem under ARUM is generally  non-concave and  would have several locally optimal solutions, even when there is only one customer zone, i.e., $|I|=1$. In contrast, we show that the cost optimization problem under the MRUM is convex and can be solved to optimality  by a nonlinear optimization solver. Furthermore, the cost optimization problem under MRUM can be formulated as a conic program, for which some off-the-shelf solvers such as CPLEX or GUROBI can efficiently handle.

\subsubsection{The Cost Optimization MCP under ARUM is Non-unimodal}
We first consider the cost optimization problem under  ARUM. The following remark states that $f^{\ARUM}(\by,\bx)$ is not unimodal in $\bx$. As a result, the MCP under ARUM, even when we only seek optimal costs, is generally not tractable to handle. 
\begin{proposition}
\textit{The objective function of the MCP under ARUM is  not unimodal, even when there is only one customer zone, i.e., $|I| = 1$. }
\end{proposition}
\proof{} We validate the remark with an example. Let us consider the following instance of $|I|=1$,
\[
f^{\ARUM}(\bx) = \frac{e^{x_1} + 1.2e^{x_2}}{1+ e^{x_1} + 1.2e^{x_2}}
\]
and a feasible set  $\cX$ defined as $\cX = \{\bx \in  [0,1]^2|\; x_1+x_2\leq 1\}$. 
Here since we fix $S$, we remove the corresponding notation for notational simplicity. Clearly, the above objective function is an instance of the MCP under ARUM with 1 customer zone and 2 locations. The MCP can be written equivalently as $$\max_{x_1 \in [0,1]} \left\{e^{x_1}+ 1.2e^{1- x_1}\right\}, $$
and it is easy to see that there are two local optimal solutions at $x_1 = 0$ and $x_1 = 1$  and the objective function achieves its maximum value at $x_1 = 0$. We illustrate the two local optimums in Fig. \ref{fig:example-fun-MCP-ARUM} below. 
Thus,  MCP has two local optimal solutions at $\bx = (1,0)$ and $\bx = (0,1)$.   
\begin{figure}[htb] 
\centering
    \includegraphics[width=0.4\linewidth]{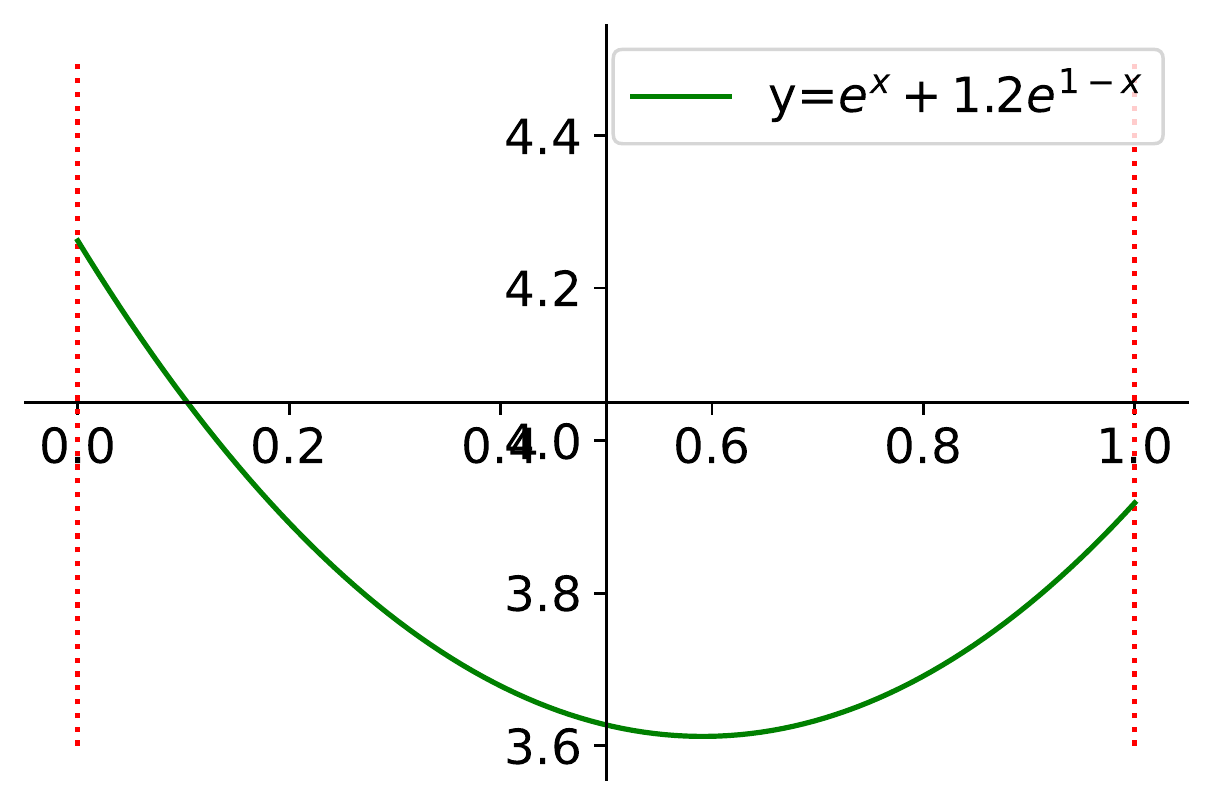} 
    \caption{Illustrative example} 
    \label{fig:example-fun-MCP-ARUM} 
\end{figure}
\endproof

In fact, when $|I|=1$, we can write the objective value of the cost optimization problem under ARUM as
\[
\max_{\bx \in \cX}\{f^{\ARUM}(\bx)\} = q_n - \frac{1}{1+ \max_{\bx\in \cX}\left\{ \sum_{i\in S} e^{a_{ni}x_{i}+b_{ni}} \right\}}
\]
where $n$ is the unique customer zone in $I$. The problem $  \max_{\bx\in \cX}\left\{ \sum_{i\in S} e^{a_{ni}x_{i}+b_{ni}} \right\}$ is a maximization problem with a convex objective function, thus several local optimums may exist.  In general, the non-unimodality of the cost optimization problem under MRUM implies
that a nonlinear optimization algorithm would end up with a sub-optimal solution. 

\subsubsection{Cost Optimization MCP under MRUM}
We consider the cost optimization problem under the MRUM framework, where the objective function can be written as 
\begin{align}
f^{\MRUM}(\bx) &=    \sum_{n\in I} q_n \frac{\sum_{i\in S} ({a_{ni} x_i +b_{ni} })}{U^c_n+ \sum_{i\in S} ({a_{ni}x_{i} + b_{ni}})}  \nonumber \\
   &=\sum_{n\in I} q_n  - \sum_{n\in I} \frac{q_n U^c_n }{U^c_n+ \sum_{i\in S} ({a_{ni}x_{i} + b_{ni}})}  \nonumber
\end{align}
It is known that a function of the form $1/(\bap^\transpose \bz)$ is convex  in $\bz$ \citep{BenaHans02,mai2020multicut}. As a result, $f^{\MRUM}(\by, \bx)$ is concave. We state this important result in Remark~\ref{prop:concave} below.
\begin{remark}[Concavity]
\label{prop:concave}
$f^{\MRUM}(\bx)$ is concave in $\bx$.
\end{remark}
The concavity implies that the cost optimization problem under MRUM can be solved to optimality by a general nonlinear optimization algorithm, e.g., gradient ascent, trust region, or linear search \citep{NoceWrig06}.

\noindent \textbf{{Conic reformulation}}.
Conic optimization refers to optimization of a linear function over conic quadratic inequalities of the form $||\bA\bz -\bb|| \leq \bc \bx + \bd$, where $||\cdot||$ is the L2 norm, and $\bA, \bb,\bc,\bd$ are matrices/vectors of appropriate sizes. Solvers for such a conic quadratic optimization problem are now available in some commercial optimization tools such as CPLEX or GUROBI.  
We will make use of the fact that such a rotated second-order cone/hyperbolic inequality of the form $z_1^2\leq z_2z_3$, for $z_1,z_2,z_3\geq 0$, can be formulated as a conic quadratic inequality  $||(2z_1,z_2-z_3)||\leq z_2+z_3$. 

In the following, we  show that the cost optimization under MRUM can be formulated as a conic quadratic optimization problem. We first let
\begin{align}
    \theta_n &= \frac{
1}{(U^c_n+\sum_{i\in S}b_{ni}) + \sum_{i\in S} a_{ni} x_{i}} \nonumber \\
w_n &= {(U^c_n+\sum_{i\in S}b_{ni}) + \sum_{i\in S} a_{ni} x_{i}} \nonumber
\end{align}
We then can formulate the cost optimization problem under MRUM as
\begin{align}
     \max_{\bx,\bw,\btheta } &\quad \sum_{n\in I} q_n - \sum_{n\in I} q_n U^c_n \theta_n \label{prob:CP-CONIC-1} \tag{\sf CP}     \\
    \mbox{subject to} &\quad  w_n = {(U^c_n+\sum_{i\in S}b_{ni}) + \sum_{i\in S} a_{ni} x_{i}} ,\;\forall n\in I \nonumber\\
    & \quad \theta_n w_n = 1 ,\;\forall n\in I \label{eq:ctr-wtheta1} \\
    &\quad \sum_{i \in S} x_i \leq C, \nonumber\\
    &\quad x_i \in [L_i,U_i], \; \forall i \in S \nonumber\\
     &\quad \bx \in \bbR^{|S|}_+,\; \bw,\btheta \in \bbR^{|I|}_+. \nonumber
 \end{align}
Constraints \eqref{eq:ctr-wtheta1} are not rotated second-order cones, but we can safely replace them by rotated second-order cones $\theta_n w_n \geq 1$. We state this result in the following proposition.
\begin{proposition}\label{prop:cp-conic}
\eqref{prob:CP-CONIC-1} can be formulated as the following conic quadratic optimization problem
\begin{align}
     \max_{\bx,\bw,\btheta } &\quad \sum_{n\in I} q_n - \sum_{n\in I} q_n U^c_n \theta_n \label{prob:CP-CONIC} \tag{\sf CP-CONIC}     \\
    \mbox{subject to} &\quad  w_n = {(1+\sum_{i\in S}b_{ni}) + \sum_{i\in S} a_{ni} x_{i}} ,\;\forall n\in I \nonumber\\
    & \quad \theta_n w_n \geq  1 ,\;\forall n\in I \nonumber \\
    &\quad \sum_{i \in S} x_i \leq C, \nonumber\\
    &\quad x_i \in [L_i,U_i], \; \forall i \in S \nonumber\\
     &\quad \bx \in \bbR^{|S|}_+,\; \bw,\btheta \in \bbR^{|I|}_+. \nonumber
\end{align}
\end{proposition}
\proof{}
We observe that  if we replace the equality constraints by inequalities $\theta_nw_n\geq 1$, 
then because $\theta_n$ and $w_n$ are non-negative and the objective function involves terms $q_n\theta_n$, the maximization problem always forces $\theta_n$ to be as small as possible. Thus, at optimum, the inequality constraints $\theta_nw_n\geq 1$ are always active, i.e.,  $\theta_nw_n= 1$.
This leads to the desired result.  
\endproof

The conic program in \eqref{prob:CP-CONIC} has $m+2|I|$ continuous variables while the original concave optimization problem  $\max_{\bx} \{ f^{\MRUM}(\bx) \}$ has only $m$ continuous variables. However, it is  expected that solving the conic program using an advanced optimization tool such as CPLEX would be faster and more scalable than directly solving the problem by a convex optimization solver, due to some advanced linear and conic programming techniques implemented in CPLEX's conic solvers \citep{bonami2015recent}. 

The cost optimization MCP problem has the form of a sum-of-ratio program. In its most general form, the problem can be written as $\max_{\bx}\{\sum_{k\in [N]} ({(\bap^k)^\transpose \bx})/((\bbt^k)^\transpose \bx) \}$, where $N\in \mathbb{N}_+$, and $\bap^k, \bbt^k$, $k\in [N]$, are vectors of the same size with $\bx$. 
Such a sum-of-ratio problem is known to be highly non-convex and challenging to be solved globally \citep{schaible1995fractional}. In our context, the cost optimization problem can be solved by convex optimization because the objective function can be written as a sum of ratios of constant numerators.   

It is also interesting to link the formulation in \eqref{prob:CP-CONIC} to the context of price optimization under discrete choice models. Under the MRUM framework, if the choice model is the MNL, the objective function of the corresponding price optimization problem contains only one ratio, thus the problem can be conveniently converted to a linear program \citep{schaible1995fractional}. If the choice model is the mixed MNL, the problem will have the form of the general sum-of-ratio program and  would be highly non-convex  and challenging to handle. 
 

\subsection{Joint Location and Cost Optimization under MRUM}
\label{sec:JOINT PROB}
We consider the joint facility location and cost optimization problem under MRUM. We will present our three approaches, two exact and one heuristic methods,  to handle the problem, namely, 
an equivalent mixed-integer conic program, a multicut outer-approximation algorithm and  a local search procedure. 


\subsubsection{{CONIC Reformulation}} 
 We will also make use of the following rotated second-order cone $x_1^2 \leq x_2x_3$, with $x_1,x_2,x_2 \geq 0$ to present conic quadratic inequality. To reformulate the nonlinear binary program \eqref{prob:MCP-mulpli} as a CONIC program, we denote
 \begin{align}
     \theta_n &=\frac{1}{U^c_n+\sum_{i\in  [m]} b_{ni} y_i + \sum_{i\in [m]} a_{ni} y_ix_{i}}\nonumber \\
     w_n &= U^c_n+\sum_{i\in  [m]} b_{ni} y_i + \sum_{i\in [m]} a_{ni} y_ix_{i}\nonumber
 \end{align}
and  write \eqref{prob:MCP-mulpli} as
\begin{align}
     \min_{\bx,\by, \bw,\btheta } &\quad \sum_{n\in I}  q_n U^c_n \theta_n \label{prob:MCP-P3}\tag{P3} \\
    \mbox{subject to} &\quad  \theta_n w_n = 1,\;\forall n\in I \nonumber\\
    & \quad w_n = U^c_n+\sum_{i\in  [m]} b_{ni} y_i + \sum_{i\in [m]} a_{ni} y_ix_{i},\;\forall n\in I\nonumber \\
    &\quad \text{Constraints \eqref{eq:ctr-capacity}-\eqref{eq:ctr-capacity-M}-\eqref{eq:ctr-x-y-j}-\eqref{eq:ctr-x-y-j-lb}}, \nonumber\\
    &\quad \bx \in \bbR^{m}_+,\; \by\in \{0,1\}^m,~ \bw,\btheta \in \bbR^{|I|}_+.   \nonumber
\end{align}
We now let  $z_{i}  = y_i x_i \; \forall i \in [m]$. Since this term involves a binary variable, the equality can be linearized using McCormick inequalities \citep{mccormick1976computability}  as
\begin{align}
    &z_i \leq U_i y_i,\; z_i \geq L_i y_i\nonumber  \\
    &z_i \leq x_i - L_i(1-y_i),\;z_i \geq x_i - U_i(1-y_i) \nonumber
\end{align}
Furthermore, since both $\theta_n$ and $w_n$ are non-negative, similarly to Proposition \ref{prop:cp-conic}, one can show that the equality $\theta_nw_n=1$ can be replaced by the inequality $\theta_nw_n\geq 1$, which is a rotated second-order cone. In summary,  we can formulate \eqref{prob:MCP-mulpli} as the following CONIC program
\begin{align}
     \min_{\bx,\by, \bz, \bw,\btheta } &\quad \sum_{n\in I}  q_n U^c_n \theta_n \label{prob:Joint-CONIC}\tag{\sf FC-CONIC} \\
    \mbox{subject to} &\quad  \theta_n w_n \geq 1,\;\forall n\in I \nonumber\\
    &\quad  w_n = 1+\sum_{i\in  [m]} b_{ni} y_i + \sum_{i\in [m]} a_{ni} z_{i},\;\forall n\in I \nonumber\\
    &\quad  \sum_{i\in [m]} z_i\leq C \label{eq:sum-z}\\
    &\quad \text{Constraints \eqref{eq:ctr-capacity-M}-\eqref{eq:ctr-x-y-j}-\eqref{eq:ctr-x-y-j-lb}}, \nonumber\\
    &\quad  z_i \leq U_i y_i,\; z_i \geq L_i y_i,\;\forall i\in [m]\label{eq:mc-1}  \\
    &\quad z_i \leq x_i - L_i(1-y_i),\;z_i \geq x_i - U_i(1-y_i),\;\forall i\in [m] \label{eq:mc-2} \\
    &\quad \bx \in \bbR^{m}_+ , \: \by \in \{0,1\}^m, \: \bz \in \bbR_+^{ m},~ \bw,\btheta \in \bbR^{|I|}_+.\nonumber
\end{align}

The mixed-integer conic program in \eqref{prob:Joint-CONIC} has $2m+2|I|$ continuous and $m$ binary variables, and has $6m + 2|I| +2$ constraints, as compared to $m$ continuous and $m$ binary variables and 2 constraints of the original joint problem \eqref{prob:MCP-mulpli}.  

\subsection{Multicut Outer-Approximation}\label{subsec:OA}
We show that an outer-approximation algorithm can be used to solve the joint facility location and cost optimization problem. We also let  $z_{i}  = y_i x_i, i \in [m]$, and linearize these terms using McCormick inequalities \citep{mccormick1976computability} to obtain the equivalent reformulation 
\begin{align}
    \max_{\bx,\by,\bz } &\quad f(\by,\bx,\bz ) =\sum_{n\in I} q_n - \sum_{n\in I} \frac{{q_nU^c_n }}{U^c_n+ \sum_{i\in [m]} {a_{ni}z_i + y_ib_{ni}}}\label{prob:MCP-mulpli-2}\tag{\sf MCP-MRUM-2}   \\
    \mbox{subject to}
    &\quad \text{Constraints \eqref{eq:ctr-capacity-M}-\eqref{eq:ctr-x-y-j}-\eqref{eq:ctr-x-y-j-lb}} \nonumber\\
    &\quad \text{Constraints \eqref{eq:sum-z}-\eqref{eq:mc-1}-\eqref{eq:mc-2}}\nonumber\\
    &\quad \bx,\bz \in \bbR^{m}_+,\: \by \in \{0,1\}^m. \nonumber
\end{align}
We can further see that the objective function of \eqref{prob:MCP-mulpli-2} is concave in ${\bz,\by}$. This remark is essential for the use of the outer-approximation method. We state it in the following remark.
\begin{remark}[Concavity]
Each component $q_n\big/({U^c_n+ \sum_{i\in [m]} {a_{ni}z_i + y_ib_{ni}}})$ of the objective function of \eqref{prob:MCP-mulpli-2} is concave in $\by$ and $\bz$, $\forall  n\in I$. 
\end{remark}
To apply the multicut outer-approximation scheme \citep{OA_Duran1986outer,mai2020multicut}, we divide the set of customer zones $I$ into
$\cL$ disjoint groups $\cD_1,\ldots,\cD_{\cL}$ such that $\bigcup_{l\in \cL}\cD_l = I$.  We then write the objective function of \eqref{prob:MCP-mulpli-2} as 
\[
f(\bx,\by,\bz) = \sum_{l\in[\cL]} \phi_l(\by, \bz)
\]
where 
\[
\phi_l(\by, \bz) = \sum_{n\in \cD_l}\left(q_n - \frac{{q_n U^c_n }}{U^c_n+ \sum_{i\in [m]} {a_{ni}z_i + y_ib_{ni}}}\right)
\]
A multicut outer-approximation algorithm will  build a piecewise linear function that outer-approximates each concave function $\phi_l(\by,\bz)$. Such a outer-approximation can be constructed by creating sub-gradient cuts  of the form $\theta_l \leq \nabla_{\by} \phi_l(\by,\bz) (\bz- \overline{\bz})  + \nabla_{\bz} \phi_l(\by,\bz) (\by- \overline{\by})   +\phi_l(\overline{\by}, \overline{\bz})$, where $(\overline{\by},\overline{\bz})$ is the current solution candidate,  and $\nabla_{\by} \phi_l(\by,\bz)$ and $\nabla_{\bz} \phi_l(\by,\bz)$ are the gradients of $\phi_l(\cdot)$ w.r.t. $\by$, $\bz$, respectively. The outer-approximation algorithm  executes an iterative procedure where at each step it adds sub-gradient cuts to the following master problem  
\begin{align}
    \underset{\bx,\by,\bz, \bw,\btheta}{\max} \qquad & \sum_{l\in [\cL]} \theta_l  & \label{prb:sub-OA}\tag{\sf sub-MOA}\\
    \text{subject to} \qquad &\quad \theta_l \leq \nabla_{\by} \phi_l(\by,\bz) (\bz- \overline{\bz}^t)  + \nabla_{\bz} \phi_l(\by,\bz) (\by- \overline{\by}^t)   +\phi_l(\overline{\by}^t, \overline{\bz}^t),~ l\in[\cL], t =1,2\ldots &\label{eq:sub-OA-cuts}\\
    &\quad \text{Constraints \eqref{eq:ctr-capacity-M}-\eqref{eq:ctr-x-y-j}-\eqref{eq:ctr-x-y-j-lb}} \nonumber\\
    &\quad \text{Constraints \eqref{eq:sum-z}-\eqref{eq:mc-1}-\eqref{eq:mc-2}}\nonumber\\
&\quad \bx,\bz \in \bbR^{m}_+,\: \by \in \{0,1\}^m,\: \btheta \in \bbR^{\cL}_+,
\end{align}
where $(\overline{\bz}^t,\overline{\by}^t)$ are the solution candidates found at the $t$-th iteration of the  outer-approximation algorithm and
the linear constraints in \eqref{eq:sub-OA-cuts} are the sub-gradient cuts added to the master problem over iterations $t=1,2,\ldots$. The algorithm keeps adding cuts and solving the master problem until its get a solution $(\bx^*,\by^*,\bz^*,\btheta^*)$ such that $\sum_{l\in [\cL]}\theta^*_l\leq \sum_{l\in [\cL]} \phi_l(\by^*,\bz^*)+\tau$, where $\tau>0$ is a stopping  threshold. The concavity of $\phi_l(\by,\bz)$, $\forall l\in[\cL]$ guarantees that  the outer-approximation algorithm will eventually converge to an optimal solution  \citep{OA_Duran1986outer,Bonami2011_BB_MIP}.
Here, we note that since the objective function is concave, one possible alternative is the outer-approximation-based Branch-and-Cut (B\&C) method \citep{Ljubic2018outer,Lin2021branch,Lin2021generalized}. We however stick to the multicut outer-approximation approach as it generally outperforms the B\&C one \citep{mai2020multicut}

\subsection{Local Search Heuristic}
\citep{Dam2021submodularity} show that a local search procedure can achieve state-of-the-art performance for the facility location MCP  under ARUM. Therefore, in this section, we explore a local search procedure to solve the joint facility and cost optimization problem under  MRUM. First, let us define
\begin{equation}
\label{prob:max-local-search}
\Phi(\by) = \max_{\bx\in \cX(\by)} \{f^{\MRUM}(\by,\bx)\},
\end{equation}
where $\cX(\by)$ is the feasible set of $\bx$ conditional on $\by$, i.e., $\cX(\by) = \{\bx\in \bbR^{m}_+|\; \sum_{i\in [m]} x_i\leq C,\; x_i\leq y_iU_i,\; x_i\geq y_i L_i,\;\forall i\in [m] \}$. According to Section \ref{sec:CP}, for a given $\by$, the value of $\Phi(\by)$ can be computed efficiently by solving a convex optimization problem, or an equivalent CONIC program. We now can  perform a local search procedure on $\Phi(\bx)$. Prior works show that the objective function of the facility location MCP is monotonic, i.e., opening new facilities always increases  the objective function. However, we show that it is not the case here in the context of the joint facility location and cost optimization problem.
\begin{proposition}
\label{prop:not-monotonic}
$\Phi(\by)$ is not necessarily monotonic increasing  in $\by$.
\end{proposition}
\proof{}
We provide a counter example to prove  the claim. Let us consider a simple example with $I = \{1\}$ (there is only one customer zone), $m=2$, $q_n=1$ for all $n\in I$, $b_{11} = b_{12} = 0$, $a_{11} = 5$, $a_{12}=1$, $U_1 = U_2=3$, $L_1 = L_2=3$, $C = 4$. Let $\by^1 = (1,0)$ and $\by^2 = (1,1)$. We see that
\begin{align}
    \Phi(\by^1) &= \max_{x_1 \in [2,3]}\left\{\frac{5x_1}{1+5x_1}\right\} = \frac{15}{16}\nonumber \\
    \Phi(\by^2) &= \max_{x_1,x_2 \in [2,3], x_1+x_2\le 4}\left\{\frac{5x_1+x_2}{1+5x_1+x_2}\right\} = \frac{12}{13}.\nonumber
\end{align}
Thus, $\Phi(\by^1) > \Phi(\by^2)$, implying that $\Phi(\by)$ is not monotonic increasing in $\by$, as desired.  
\endproof

Intuitively, when a new facility at location $i$ is built and operated, it would require a minimum cost  $L_i>0$, which would reduce the investment/operational  costs for other  facilities. The objective  function then would decrease if the new  facility is much less attractive to the customer than the other ones. On the other hand, if the total budget $C$ is sufficiently large, e.g., $C\geq \sum_{i\in[m]} U_i$, then we can see that $\Phi(\by)$ is monotonic increasing in $\by$. 

Proposition \ref{prop:not-monotonic} implies that an optimal solution to \eqref{prob:max-local-search} would not necessarily reach the maximum capacity $K$. Moreover, in prior MCP work, \cite{Dam2021submodularity} show that a simple greedy local search based on adding locations to an empty set of locations can   guarantee $(1-1/e)$ approximation solutions. 
This would be not the case for  \eqref{prob:max-local-search} due to the non-monotonicity \citep{Nemhauser1978analysis}.

\cite{Dam2021submodularity} propose a three-step local search procedure that performs very well for the ARUM-based facility location MCP. The second step of this algorithm requires gradient information to direct the local search. To adapt this in the context of the joint facility location and cost optimization problem, we  need to compute the gradients of $\Phi(\by)$ w.r.t.  any $\by$. To this end, let us consider a Lagrange dual of the maximization problem in \eqref{prob:max-local-search}.
\begin{equation}
\label{prob:Lagrange-dual}
L(\bx,\lambda,\bgamma^U,\bgamma^L |\by) = f^{\MRUM}(\by,\bx) - \lambda \left(\sum_{i\in[m]} x_i-C\right) - \sum_{i\in [m]}\gamma^U_i(x_i-y_iU_i) + \sum_{i\in [m]}\gamma^L_i(x_i-y_iL_i).
\end{equation}
For a given $\by\in\cY$ with $\cY$ is the feasible set of $\by$, i.e., $\cY = \{ \by \in \{ 0,1 \}^m|\sum_{i \in [m]} y_i \leq K \}$, let $\bx^*(\by), \lambda^*(\by), \gamma^{U*}_i(\by)$ and $\gamma^{L*}_i(\by)$, $\forall i\in [m]$ be the saddle point of $L(\bx,\lambda,\bgamma^U,\bgamma^L |\by)$.
The following proposition gives a formula for the first-order derivatives of $\Phi(\by)$ w.r.t. $\by$.
\begin{proposition}
\label{prop:derivative-phi}
The first-order derivatives of $\Phi(\by)$ is given as 
\[
\frac{\partial \Phi(\by)}{\partial y_i} =\frac{\partial f^{\MRUM}(\by,\bx^*) }{\partial y_i} + \gamma^{U*}_i(\by) U_i - \gamma^{L*}_i(\by) L_i,\;\forall i\in [m].  
\]
\end{proposition}

Proposition \ref{prop:derivative-phi} tells us that we need to compute the Lagrange multipliers $\lambda^*, \bgamma^{U*}$ and $\bgamma^{L*}$ to compute the derivatives of $\Phi(\by)$ w.r.t. $\by$. This can be done by solving the Lagrange duality $\min_{\lambda,\bgamma^{U},\bgamma^{L}} \max_{\bx} L(\bx,\lambda,\bgamma^U,\bgamma^L |\by)$, which would be expensive. Alternatively, we can  directly solve the convex problem $\max_{\bx\in \cX(\by)} f^{\MRUM}(\by,\bx)$ by an existing optimization solver and infer the Lagrange multipliers  from the  optimal solution obtained. That is, let $\bx^*$ be the unique optimal solution to $\max_{\bx\in \cX(\by)} f^{\MRUM}(\by,\bx)$. From \eqref{eq:dual-2} we have 
\[
\frac{\partial f^{\MRUM}(\by,\bx^*)}{\partial x_i} -\lambda^* - \gamma^{U*}_i + \gamma^{L*}_i = 0,\;\forall i\in[m].
\]
To infer $\lambda^*,\bgamma^{U*}$ and $\bgamma^{L*}$,  we consider the following  cases:
\begin{itemize}
    \item  If $\sum_{i\in[m]} x^*_i <C$, then from \eqref{eq:dual-3} we see that $\lambda^* = 0$. Then using the complementary slackness conditions in \eqref{eq:dual-4}-\eqref{eq:dual-5}  we can infer the values of $\gamma^{U*}_i$ and $\gamma^{L*}_i$ as follows. For each $x_i^*$ such that $y_i>0$ we have
\[
\begin{cases}
\text{If } x^*_i = y_i U_i\text{ then }\gamma^{L*}_i = 0  \text{ and }\gamma^{U*}_i = \frac{\partial f^{\MRUM}(\by,\bx^*)}{\partial x_i}\\
\text{If } x^*_i = y_iL_i\text{ then }\gamma^{U*}_i = 0  \text{ and }\gamma^{L*}_i = -\frac{\partial f^{\MRUM}(\by,\bx^*)}{\partial x_i}\\
\text{If } y_iL_i< x^*_i< y_iU_i\text{ then }\gamma^{U*}_i =\gamma^{L*}_i=0.
\end{cases}
\]
\item If $\sum_{j}x^*_i = C$ and  if there is $i^*\in [m]$ such that $y_{i^*}L_{i^*}< x^*_{i^*} < y_{i^*}U_{i^*}$, then $\gamma^{U*}_{i^*} = \gamma^{L*}_{i^*} = 0$. Thus, $$\lambda^* =\frac{\partial f^{\MRUM}(\by,\bx^*)}{\partial x_{i^*}} $$
and for any $i\in [m]$ we have
\[
\gamma^{U*}_i - \gamma^{L*}_i  = \frac{\partial f^{\MRUM}(\by,\bx^*)}{\partial x_i} -\lambda^* \myeq{} t_i
\]
and 
\[
\begin{cases}
\text{If } x^*_i = y_i U_i\text{ then }\gamma^{L*}_i= 0  \text{ and }\gamma^{U*}_i = t_i\\
\text{If } x^*_i = y_iL_i\text{ then }\gamma^{U*}_i = 0  \text{ and }\gamma^{L*}_i = -t_i\\
\text{If } y_iL_i< x^*_i < y_iU_i\text{ then }\gamma^{U*}_i =\gamma^{L*}_i=0.
\end{cases}
\]
\item If $\sum_{i\in[m]}x^*_i = C$,  and $x^*_i=y_iU_i$ or $x^*_i=y_iL_i$ for all $i\in[m]$. This means that $C$ can be written as a sum of  some lower bounds $L_i$ and some upper bounds $U_i$. In this situation,  the Lagrange multipliers cannot be uniquely identified using the KKT conditions.  We can simply add some small amount to $C$ such that it no-longer can be written as  a sum of some upper and lower bounds $L_i, U_i$. This allows us to switch the situation to the two above cases, for which the Lagrange multipliers can be well identified. Note that, by adding  a small amount to $C$, we slightly modify the objective function $\Phi(\by)$, but we expect that if the amount added to $C$ is relatively small, we still can obtain a good approximation of the gradients of $\Phi(\by)$. 
\end{itemize}

The local search algorithm consists of three steps. In the first step, we start with an empty solution and keep adding locations to the current solution, taking locations that increase the objective function $\Phi(\by)$ the most. The first step stops  when the solution candidate reaches  the capacity $K$ or it cannot find a better solution. The second step is a gradient-based local search. The gradient information given in Proposition \ref{prop:derivative-phi} is made use to approximate the nonlinear objective function $f^{\MRUM}$ and find the next solution candidates.  The last step is a simple exchanging procedure trying to swap a location in the current choice set with some other locations outside the current choice set of locations to improve the solution found after the second step. We refer the reader to \cite{Dam2021submodularity} for more details. 

The above local search procedure runs fast and returns good solutions for previous competitive facility location problems \citep{Dam2021submodularity,thuy2021robust}, but it would be not the case in the context of the joint problem because each iteration of the local search requires to compute several values of $\Phi(\by)$, thus requires to solve several cost optimization problems. In other words, the local search algorithm developed for the joint location and cost optimization under MRUM  would be expensive and not as competitive as the other exact methods discussed above (i.e., the CONIC reformulation and the outer-approximation algorithm). We will show this clearly in the experimental section below.

\section{Numerical Experiments}
\label{sec:experiments}
In this section, we provide experimental results for the two problems we study, i.e. the cost optimization MCP and the joint location and cost optimization MCP problems under MRUM. We first present our experimental settings and then provide our comparison results. 
\subsection{Experimental Settings}
We randomly generated a dataset for each problem to evaluate the performance of the  algorithms proposed. The dataset for the cost optimization problem has totally 240 instances, in which we generate 24 instances for each pair $(|I|,m)$ in the set \{(100, 100), (100, 1000), (100, 3000), (200, 2000), (400, 1000), (800, 400), (800, 800), (1000, 1000), (2000, 2000), (5000, 1000) \}, and the dataset for the joint location  and cost optimization problem has 216 instances of the same values of $(|I|,m)$ as in the cost optimization problem, but excludes the pair (2000, 2000)  since instances of this size become too expensive to be solved by any approaches within the time budget. 


We set the maximum budget to spend on the new facilities  as $C \in \{0.2m, 0.5m, 0.7m\}$ and the maximum number of facilities that can be opened as $M \in \{ 0.2m, 0.5m, 0.7m \}$, recalling that $m$ is the total number of possible locations. As describe above, given $n\in I$ and $i\in [m]$, the parameter $a_{ni}$ represents  the sensitivity of customer zone $n$ w.r.t the cost spent on facility $i$, and the parameter $b_{ni}$ represents other factors that affect customers' choice decisions, e.g. distances to the locations. We randomly and uniformly generate $a_{ni}$ from [0.5,1.5]. For $b_{ni}$, we also randomly generate values from [1,10].
The parameter $q_n$ representing the number of customers in customer zone $n$ is also generated randomly 
and  uniformly from two intervals  $[1,10]$ and $[90,100]$. The two parameters $L_i$ and $U_i$, $\forall i\in [m]$, (i.e., the lower and upper bounds for the cost to spend on facility $i\in [m]$) are also generated uniformly from  two closed intervals satisfying that $\sum_{i\in[m]}L_i \leq C$ for some sets of $i\in[m]$. We note that the dataset for the cost optimization MCP problem needs to satisfy $\sum_{i\in[m]}L_i \leq C$ to ensure that the problem always has a feasible solution.

The experiments are conducted on a PC with processors AMD Ryzen 3-3100 CPU @ 3.60 GHz, RAM of 16 gigabytes, and operating system Window 10. We use MATLAB 2020 to implement and run the convex optimization algorithms for the cost optimization problem. We use C++ and link to IBM ILOG-CPLEX 12.10 (under default settings) to solve CONIC programs and the master problem of the multicut outer-approximation algorithm. The CPU time limit for each instance is 600 seconds, i.e., we stop the algorithms if they exceed the time budget and report the best solutions found.
\subsection{Cost Optimization MCP}
\label{subsec:results-cost}
In this section, we provide comparison results for the cost optimization MCP. We solve the instances generated above by the two approaches discussed in Section \ref{sec:CP}, i.e.,  using a general convex optimization solver, and solving the equivalent CONIC program by CPLEX.
We denote the two approaches as CONIC  and CONVEX (i.e., solving the convex optimization problem by  \textit{fmincon}).  
Table \ref{table1} below reports the numbers of instances solved with the best objectives and the average CPU times. Each row of this table reports  averages over 8 instances grouped by$(|I|, m)$, and $C$. We highlight  in bold the largest numbers of instances solved with best objective values, noting that both CONIC and CONVEX are exact methods, thus if an instance is completed by any of the two approaches within the time budget, we know that the instance is solved to optimality.

\begin{table}[htb]\footnotesize
\centering
\caption{Comparison results for the cost optimization MCP  instances, 8 instances per row.}
\label{table1}
\begin{tabular}{lllll|ll} 
\multirow{2}{*}{$|I|$} & \multirow{2}{*}{$m$} & \multirow{2}{*}{$C$} & \multicolumn{2}{c|}{\begin{tabular}[c]{@{}c@{}}\# Instances \\with best objective\end{tabular}} & \multicolumn{2}{c}{\begin{tabular}[c]{@{}c@{}}Average \\CPU time (s)\end{tabular}}  \\ 
\cline{4-7}
    &   &   & CONIC& CONVEX     & CONIC  & CONVEX     \\ 
\hline
100 & 100& 20& \textbf{8 } & \textbf{8 }& 0.05   & 0.11\\ 
\hline
100 & 100& 50& \textbf{8 } & \textbf{8 }& 0.05   & 0.09\\ 
\hline
100 & 100& 70& \textbf{8 } & \textbf{8 }& 0.05   & 0.10\\ 
\hline
100 & 1000& 200& \textbf{8 } & \textbf{8 }& 0.41   & 6.10\\ 
\hline
100 & 1000& 500& \textbf{8 } & \textbf{8 }& 0.42   & 7.36\\ 
\hline
100 & 1000& 700& \textbf{8 } & \textbf{8 }& 0.43   & 6.96\\ 
\hline
800 & 400& 80& \textbf{8 } & \textbf{8 }& 1.84   & 5.48\\ 
\hline
800 & 400& 200& \textbf{8 } & \textbf{8 }& 1.85   & 6.10\\ 
\hline
800 & 400& 280& \textbf{8 } & \textbf{8 }& 1.79   & 6.32\\ 
\hline
800 & 800& 160& \textbf{8 } & \textbf{8 }& 3.77   & 50.48      \\ 
\hline
800 & 800& 400& \textbf{8 } & \textbf{8 }& 3.76   & 52.97      \\ 
\hline
800 & 800& 560& \textbf{8 } & \textbf{8 }& 3.77   & 57.25      \\ 
\hline
400 & 1000& 200& \textbf{8 } & \textbf{8 }& 1.55   & 22.72      \\ 
\hline
400 & 1000& 500& \textbf{8 } & \textbf{8 }& 1.53   & 23.69      \\ 
\hline
400 & 1000& 700& \textbf{8 } & \textbf{8 }& 1.56   & 24.31      \\ 
\hline
1000& 1000& 200& \textbf{8 } & \textbf{8 }& 7.30   & 112.53     \\ 
\hline
1000& 1000& 500& \textbf{8 } & \textbf{8 }& 7.18   & 121.48     \\ 
\hline
1000& 1000& 700& \textbf{8 } & \textbf{8 }& 7.44   & 129.71     \\ 
\hline
200 & 2000& 400& \textbf{8 } & \textbf{8 }& 1.85   & 55.35      \\ 
\hline
200 & 2000& 1000& \textbf{8 } & \textbf{8 }& 1.89   & 55.54      \\ 
\hline
200 & 2000& 1400& \textbf{8 } & \textbf{8 }& 1.85   & 50.73      \\ 
\hline
100 & 3000& 600& \textbf{8 } & \textbf{8 }& 1.84   & 60.93      \\ 
\hline
100 & 3000& 1500& \textbf{8 } & \textbf{8 }& 1.88   & 63.72      \\ 
\hline
100 & 3000& 2100& \textbf{8 } & \textbf{8 }& 1.90   & 49.79      \\ 
\hline
5000& 1000& 200& \textbf{8}  & 4     & 136.08 & 621.93     \\ 
\hline
5000& 1000& 500& \textbf{8}  & 2     & 116.72 & 670.42     \\ 
\hline
5000& 1000& 700& \textbf{8}  & 0     & 114.81 & 661.56     \\ 
\hline
2000& 2000& 400& \textbf{8 } & 0     & 48.40  & 600.00     \\ 
\hline
2000& 2000& 1000& \textbf{8 } & 0     & 46.44  & 600.00     \\ 
\hline
2000& 2000& 1400& \textbf{8 } & 0     & 47.55  & 600.00     \\ 
\hline\hline
\multicolumn{3}{c}{Average}    & \textbf{8}  & 7.1   & &\\
\hline
\end{tabular}
\end{table}

Table \ref{table1} shows that the CONIC approach outperforms the general convex solver in terms of both number of instances with the best objective values and average CPU time. The CONIC approach gives the best objective values for all the instances in
all the 30 groups of instances while the CONVEX  just gives the largest numbers of best objective values in 24/30 groups of instances. On the other hand, the average CPU times required by CONIC are significantly less than the time taken by  CONVEX. 
{For some medium and large instances (i.e., when $m>400$), the CONIC approach is about 5 to 30 times faster than the CONVEX. 
More specifically, 
for instances of $(|I|,m) = (100,1000)$, the average CPU times required by  CONVEX are about 15 times larger than those required by the CONIC. CONIC is even about 30 times faster than CONVEX for instances of  $m = 3000$.}
In addition, we can see that an increase  in the  number of zones $|I|$ or the number of locations $m$ produces remarkable changes in performance for both CONIC and CONVEX, while a change in $C$ has almost no effect.

\subsection{Joint Location  and Cost Optimization MCP} 
We now report numerical results for the joint location and cost optimization problem.
We perform a comparison between the three approaches proposed above, namely, the CONIC reformulation solved by the CPLEX's conic solver \citep{bonami2015recent}, the multicut outer-approximation algorithm (denoted as MOA), and the local search heuristic (denoted as LS).
We note that one of the most important parameters affecting the performance of the MOA algorithm is the number of groups $T$ (i.e, the number of cuts  added to the master problem at each iteration). Therefore, we test our MOA algorithm with $T \in \{ 1, 5\}$. A more detailed analysis on the performance of MOA when $T$ varies will be provided later. For the LS approach,  the computation of each $\Phi(\by)$ requires to solve a cost optimization problem over $\bx$. To this end, we use the CONIC approach (i.e. the better one among the two methods proposed to solve the cost optimization problem - see Section \ref{subsec:results-cost} for instance).

We first report the comparison results  for our small-sized and medium-sized instances in Table \ref{table2}, where each row of the table reports the averages for 8 instances 
grouped by $|I|,m, C,$ and $K$. 
More results with large-sized instances can be found in the appendix.
There are 45 groups of instances in total. 
 Here, we know that the CONIC and MOA are global approaches, thus if they stop before exceeding the time budget (600 seconds), then we report that the corresponding instances  are solved to optimality. In contrast, since LS is heuristic, we do not include this approach in the ``\textit{solved to optimality}'' columns. We also take the best solutions found by the three approaches (CONIC, MOA, and LS) and report the numbers of instances solved with the best objectives. 
We also indicate the largest numbers of instances solved to optimality and the largest numbers of instances solved  with the best objective values in bold.

Table \ref{table2} shows that the MOA algorithms seem to achieve the best performance for all the instances considered in terms of the number of instances solved to optimality;
they provide optimal solutions for all the instances of each group, while the CONIC approach only gives the optimal solutions for 22/45 {groups of instances}. We also observe that the CONIC approach is less efficient in cases where $K = m/5$. {In particular, for larger problem instances (e.g., $|I| = 800$), on average, the CONIC approach only solves less than 3/8 instances to optimality within 600 seconds.}
In overall, the CONIC approach can only  solve about 79.6\% instances  to optimality while the MOA algorithms give optimal solutions for 100\% instances. However, when comparing the number of instances with the best objective values, we observe that, on average, both CONIC and MOA are able to return the best values for all the instances. This is because CONIC exceeds the time budget and cannot confirm that the solutions found are optimal. The LS algorithm is much worse, as compared to the other exact approaches, as it only gives the best solutions to only 5.6\% of the instances.

\begin{table}\scriptsize
\centering
\caption{Comparison results for the joint location and cost optimization problem, for small-sized and medium-sized instances grouped by $(m,|I|,C,K)$ (8 instances per row); CN stands for the CONIC approach and M1, M5 stand for the MOA algorithms with $T=1,5$, respectively.}
\label{table2}
\resizebox{\linewidth}{!}{%
\begin{tabular}{lllllll|lll|lll} 
\hline
\multirow{2}{*}{$m$} & \multirow{2}{*}{$|I|$} & \multirow{2}{*}{$C$} & \multirow{2}{*}{$K$} & \multicolumn{3}{c|}{\begin{tabular}[c]{@{}c@{}}\# Instances solved\\~to optimality\end{tabular}} & \multicolumn{3}{c|}{\begin{tabular}[c]{@{}c@{}}\# Instances with\\best objectives\end{tabular}} & \multicolumn{3}{c}{Average CPU time (s)} \\ 
\cline{5-13}
 &  &  &  & CN & M1 & M5 & CN & M1/M5 & LS & CN & M1 & M5 \\ 
\hline
100 & 100 & 20 & 20 & \textbf{8 } & \textbf{8 } & \textbf{8 } & \textbf{8} & \textbf{8} & \textbf{8} & 1.0 & \textbf{0.3} & \textbf{0.3} \\ 
\hline
100 & 100 & 20 & 50 & \textbf{8 } & \textbf{8 } & \textbf{8 } & \textbf{8} & \textbf{8} & \textbf{8} & 1.2 & \textbf{0.1} & \textbf{0.1} \\ 
\hline
100 & 100 & 20 & 70 & \textbf{8 } & \textbf{8 } & \textbf{8 } & \textbf{8} & \textbf{8} & 7 & 1.1 & \textbf{0.0} & 0.1 \\ 
\hline
100 & 100 & 50 & 20 & 7 & \textbf{8 } & \textbf{8 } & \textbf{8} & \textbf{8} & \textbf{8} & 77.2 & \textbf{1.1} & \textbf{1.1} \\ 
\hline
100 & 100 & 50 & 50 & \textbf{8 } & \textbf{8 } & \textbf{8 } & \textbf{8} & \textbf{8} & \textbf{8} & 2.3 & \textbf{0.3} & 0.4 \\ 
\hline
100 & 100 & 50 & 70 & \textbf{8 } & \textbf{8 } & \textbf{8 } & \textbf{8} & \textbf{8} & \textbf{8} & 3.2 & \textbf{0.1} & \textbf{0.1} \\ 
\hline
100 & 100 & 70 & 20 & \textbf{8 } & \textbf{8 } & \textbf{8 } & \textbf{8} & \textbf{8} & \textbf{8} & 7.1 & 1.3 & \textbf{1.0} \\ 
\hline
100 & 100 & 70 & 50 & \textbf{8 } & \textbf{8 } & \textbf{8 } & \textbf{8} & \textbf{8} & \textbf{8} & 1.1 & \textbf{0.2} & 0.3 \\ 
\hline
100 & 100 & 70 & 70 & \textbf{8 } & \textbf{8 } & \textbf{8 } & \textbf{8} & \textbf{8} & \textbf{8} & 1.1 & \textbf{0.0} & 0.1 \\ 
\hline
400 & 800 & 80 & 80 & 1 & \textbf{8 } & \textbf{8 } & \textbf{8} & \textbf{8} & 0 & 571.7 & \textbf{32.2} & 33.1 \\ 
\hline
400 & 800 & 80 & 200 & 7 & \textbf{8 } & \textbf{8 } & \textbf{8} & \textbf{8} & 0 & 133.2 & \textbf{2.1} & 3.6 \\ 
\hline
400 & 800 & 80 & 280 & \textbf{8 } & \textbf{8 } & \textbf{8 } & \textbf{8} & \textbf{8} & 0 & 32.5 & \textbf{0.3} & 0.4 \\ 
\hline
400 & 800 & 200 & 80 & 1 & \textbf{8 } & \textbf{8 } & \textbf{8} & \textbf{8} & 4 & 530.2 & \textbf{84.0} & 93.6 \\ 
\hline
400 & 800 & 200 & 200 & \textbf{8 } & \textbf{8 } & \textbf{8 } & \textbf{8} & \textbf{8} & 4 & 33.5 & \textbf{5.1} & 9.9 \\ 
\hline
400 & 800 & 200 & 280 & 7 & \textbf{8 } & \textbf{8 } & \textbf{8} & \textbf{8} & 4 & 105.0 & \textbf{0.4} & 0.5 \\ 
\hline
400 & 800 & 280 & 80 & 1 & \textbf{8 } & \textbf{8 } & \textbf{8} & \textbf{8} & 4 & 600.0 & \textbf{92.1} & 104.1 \\ 
\hline
400 & 800 & 280 & 200 & \textbf{8 } & \textbf{8 } & \textbf{8 } & \textbf{8} & \textbf{8} & 4 & 35.4 & \textbf{7.8} & 18.0 \\ 
\hline
400 & 800 & 280 & 280 & 7 & \textbf{8 } & \textbf{8 } & \textbf{8} & \textbf{8} & 4 & 110.9 & \textbf{0.7} & 4.2 \\ 
\hline
800 & 800 & 160 & 160 & 0 & \textbf{8 } & \textbf{8 } & \textbf{8} & \textbf{8} & 1 & 600.0 & \textbf{92.5} & 124.0 \\ 
\hline
800 & 800 & 160 & 400 & 5 & \textbf{8 } & \textbf{8 } & \textbf{8} & \textbf{8} & 0 & 274.6 & \textbf{6.5} & 9.7 \\ 
\hline
800 & 800 & 160 & 560 & \textbf{8 } & \textbf{8 } & \textbf{8 } & \textbf{8} & \textbf{8} & 0 & 104.2 & 2.4 & \textbf{2.2} \\ 
\hline
800 & 800 & 400 & 160 & 2 & \textbf{8 } & \textbf{8 } & \textbf{8} & \textbf{8} & 4 & 474.3 & 315.5 & \textbf{213.8} \\ 
\hline
800 & 800 & 400 & 400 & 6 & \textbf{8 } & \textbf{8 } & \textbf{8} & \textbf{8} & 4 & 233.0 & \textbf{16.1} & 43.4 \\ 
\hline
800 & 800 & 400 & 560 & 7 & \textbf{8 } & \textbf{8 } & \textbf{8} & \textbf{8} & 4 & 160.0 & \textbf{3.0} & 7.0 \\ 
\hline
800 & 800 & 560 & 160 & 2 & \textbf{8 } & \textbf{8 } & \textbf{8} & \textbf{8} & \textbf{8} & 460.8 & 333.2 & \textbf{275.5} \\ 
\hline
800 & 800 & 560 & 400 & \textbf{8 } & \textbf{8 } & \textbf{8 } & \textbf{8} & \textbf{8} & 4 & 101.2 & \textbf{15.7} & 38.3 \\ 
\hline
800 & 800 & 560 & 560 & \textbf{8} & \textbf{8} & \textbf{8} & \textbf{8} & \textbf{8} & 4 & 76.9 & \textbf{4.3} & 6.3 \\ 
\hline
1000 & 100 & 200 & 200 & 4 & \textbf{8} & \textbf{8} & \textbf{8} & \textbf{8} & 4 & 303.9 & \textbf{7.0} & 10.6 \\ 
\hline
1000 & 100 & 200 & 500 & \textbf{8 } & \textbf{8 } & \textbf{8 } & \textbf{8} & \textbf{8} & 2 & 12.0 & \textbf{1.2} & 3.1 \\ 
\hline
1000 & 100 & 200 & 700 & \textbf{8 } & \textbf{8 } & \textbf{8 } & \textbf{8} & \textbf{8} & 0 & 3.6 & \textbf{0.9} & 2.5 \\ 
\hline
1000 & 100 & 500 & 200 & 7 & \textbf{8 } & \textbf{8 } & \textbf{8} & \textbf{8} & \textbf{8} & 97.3 & \textbf{13.0} & 18.8 \\ 
\hline
1000 & 100 & 500 & 500 & \textbf{8 } & \textbf{8 } & \textbf{8 } & \textbf{8} & \textbf{8} & 4 & 4.5 & \textbf{2.5} & 4.4 \\ 
\hline
1000 & 100 & 500 & 700 & \textbf{8 } & \textbf{8 } & \textbf{8 } & \textbf{8} & \textbf{8} & 4 & 3.9 & \textbf{1.2} & 3.9 \\ 
\hline
1000 & 100 & 700 & 200 & \textbf{8 } & \textbf{8 } & \textbf{8 } & \textbf{8} & \textbf{8} & \textbf{8} & \textbf{7.9} & 20.4 & 20.0 \\ 
\hline
1000 & 100 & 700 & 500 & \textbf{8 } & \textbf{8 } & \textbf{8 } & \textbf{8} & \textbf{8} & \textbf{8} & 4.0 & \textbf{2.1} & 5.1 \\ 
\hline
1000 & 100 & 700 & 700 & \textbf{8 } & \textbf{8 } & \textbf{8 } & \textbf{8} & \textbf{8} & 4 & 3.6 & 1.4 & \textbf{1.3} \\ 
\hline
1000 & 400 & 200 & 200 & 1 & \textbf{8 } & \textbf{8 } & \textbf{8} & \textbf{8} & 4 & 532.2 & \textbf{35.8} & 56.4 \\ 
\hline
1000 & 400 & 200 & 500 & 7 & \textbf{8 } & \textbf{8 } & \textbf{8} & \textbf{8} & 0 & 151.3 & \textbf{2.9} & 5.8 \\ 
\hline
1000 & 400 & 200 & 700 & 7 & \textbf{8 } & \textbf{8 } & \textbf{8} & \textbf{8} & 0 & 101.6 & \textbf{1.6} & 1.7 \\ 
\hline
1000 & 400 & 500 & 200 & 7 & \textbf{8 } & \textbf{8 } & \textbf{8} & \textbf{8} & \textbf{8} & 175.8 & 100.2 & \textbf{84.4} \\ 
\hline
1000 & 400 & 500 & 500 & 7 & \textbf{8 } & \textbf{8 } & \textbf{8} & \textbf{8} & 4 & 102.4 & \textbf{6.0} & 14.4 \\ 
\hline
1000 & 400 & 500 & 700 & \textbf{8 } & \textbf{8 } & \textbf{8 } & \textbf{8} & \textbf{8} & 4 & 41.6 & 3.8 & \textbf{2.6} \\ 
\hline
1000 & 400 & 700 & 200 & 4 & \textbf{8 } & \textbf{8 } & \textbf{8} & \textbf{8} & \textbf{8} & 332.9 & 179.3 & \textbf{148.9} \\ 
\hline
1000 & 400 & 700 & 500 & 7 & \textbf{8 } & \textbf{8 } & \textbf{8} & \textbf{8} & 4 & 108.7 & \textbf{6.2} & 14.0 \\ 
\hline
1000 & 400 & 700 & 700 & 7 & \textbf{8 } & \textbf{8 } & \textbf{8} & \textbf{8} & 4 & 101.0 & 11.7 & \textbf{5.4} \\ 
\hline\hline
\multicolumn{4}{l}{Average} & 6.37 & 8 & 8 & 8 & 8 & 4.49 & &  &  \\
\hline
\end{tabular}
}
\end{table}

We now turn our attention to the average CPU times. In Table \ref{table2} we indicate in bold the smallest average CPU times. We first remark that 
 LS has the worst performance as it  always requires more than 600 seconds to solve 43/45 groups of instances.
This is contrary to what has been observed for other MCP facility location problems \citep{Dam2021submodularity} but is to be expected, as we discuss earlier, each step of the LS requires to solve numerous cost optimization problems, thus is
 much more expensive than when only the selection of locations is considered. 
Moreover, we also observe that the  MOA algorithms  outperform the other approaches (CONIC and LS) in terms of average CPU time for instances of $(|I|,m) \in$ \{(100,100), (800,400), (800,800),(400,1000)\} and they are just slightly slower than CONIC for {6 instances} of $(|I|,m) = (100,1000)$. If we count the numbers of groups for which each approach requires the shortest average CPU times (number of values in bold on each column), we can see that M1 (i.e. single-cut outer-approximation)  performs the best with 35/45 times dominating the other approaches, followed  M5 (13/45 times). In particular, the CONIC  only performs the best for only 1/45 group of instances and LS exceeds the time budget for most the cases and is always the slowest approach.

\begin{figure}[H] 
\centering
    \includegraphics[scale=0.6]{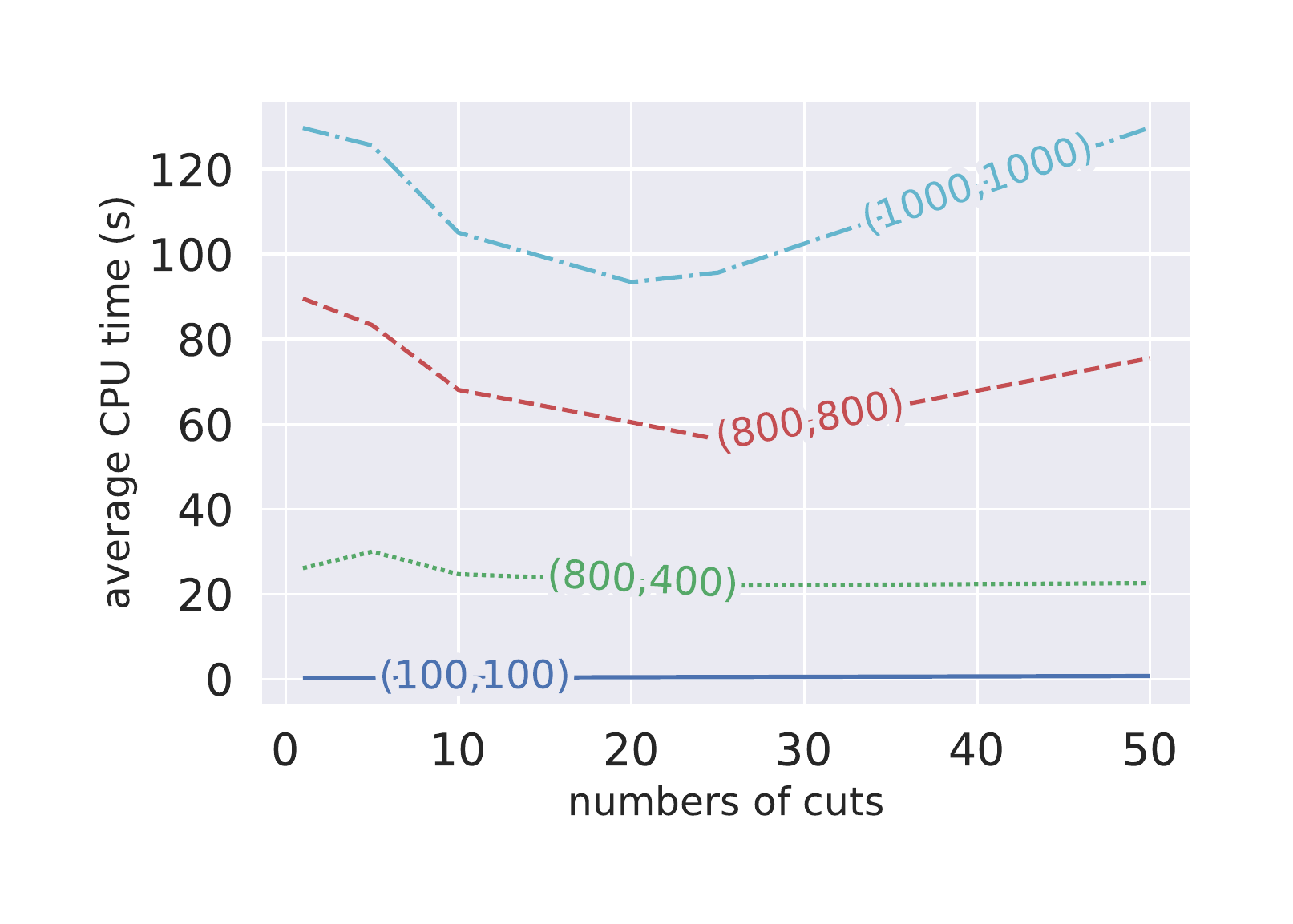}
    \caption{Average CPU times (grouped by $|I|$ and $m$) of the MOA as functions of the number of cuts per iteration $T$.}
    \label{fig:compareOutersize} 
\end{figure}

\begin{figure}[H] 
\centering
    \includegraphics[scale=0.38]{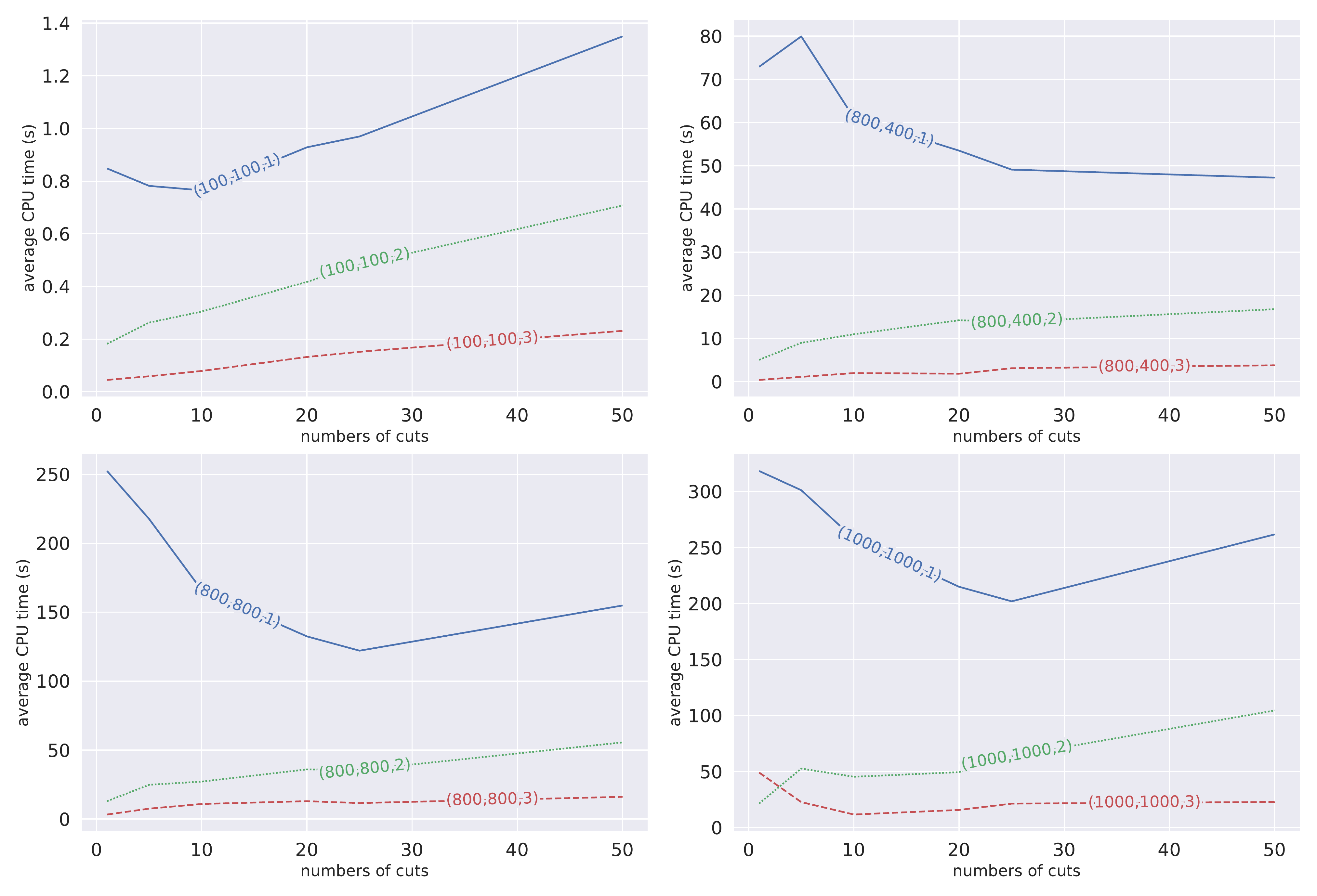}
    \caption{Average CPU times (grouped by $|I|,m$  and $K$) of the MOA algorithm as functions of the number of cuts per iteration $T$.}
    \label{compareM}
\end{figure}

Since  MOA algorithms perform the best for the joint location and cost optimization instances,  and their performance would depend on the number of cuts added to the master problem at each iteration (i.e., $T$), we provide additional results to analyse how $T$ affects the running time of the MOA. 
 To this end, we let $T$ vary in  $ \{1,5,10,20, 25,50\}$ and plot the corresponding average running times in Figure \ref{fig:compareOutersize} for  four groups of instance of sizes $(|I|,m) \in$  \{(100,100), (800,400), (800,800), (1000,1000)\}, 
 where each line represents the average CPU time of 432 instances of the same $(|I|,m)$.
 We also plot in Figure \ref{compareM} the average CPU times  for each group $(|I|,m)$ with different values of $K$, i.e., $K \in \{0.2m, 0.5m,0.7m\}$. Each line on each sub-figure  corresponds to 144 instances.
We see from the Figure \ref{fig:compareOutersize} that the MOA algorithm with $T=25$ seems to have the lowest average running times. Moreover, Figure \ref{compareM} shows that the MOA algorithms with small $T$ (i.e, $T=1$ or $T = 5$) seem to achieve good performance 
 for $K \in \{ 0.5m,0.7m\}$. However, the algorithm becomes slower for instances of $K=0.2m$. Such instances seem to be better solved by MOA with $T$ from 10 to 25. 
We also observe that, when $T$ becomes too large (e.g. $T\ge 50$), the MOA algorithm becomes slow and often encounters  out-of-memory issues.
This is to be expected as if $T$ is large, the master problem grows quickly in size, and will become much more expensive to solve. 
All these  are consistent with observations in \citep{mai2020multicut}. It is interesting to note that the MOA algorithm with $T=25$ seems  to performs  the best in terms of CPU time, but the results in  Table  \ref{table2} 
show that the MOA with $T=1$ seems to be more robust, in the sense that it gives the highest number of instances solved to optimality. This can be explained by the fact that MOA with a high value of $T$ would better explore the structure of the problem and require less iterations to converge. On the contrary, MOA with $T=1$ will require more iterations, but the cost to perform one iteration will be cheaper. Thus, the MOA with high $T$ would finish  small- or medium-size instances faster, but for large-size instances (the results are reported in the appendix), the master problem is getting expensive quickly and  the MOA would not be  able to terminate the instances within the time budget. 

In summary, for the cost optimization problem, we observe that our proposed methods (CONIC and CONVEX) are able to solve large-scale instances to optimality in fairly short CPU times, and
the CONIC approach performs better than  CONVEX. For the joint location and cost optimization problem, the MOA with $T=1$ performs the best in terms of number of instances solved to optimality.
In particular, the MOA algorithms  can efficiently handle instances of thousands of locations and customer zones. 
The LS approach, even-though performs well in prior MCP work, performs the worst.

\section{Conclusion}
\label{sec:concl}
In this paper, we studied a joint facility location and cost optimization problem under RUM models. We examined the two popular RUM frameworks in  the discrete choice literature, namely the ARUM and MRUM. 
We showed that the two frameworks can well approximate each other on the context of MCP cost optimization.
In the optimization perspective, we showed that, while the cost optimization under ARUM is highly non-convex  and would be challenging to handle,  the cost optimization problem under MRUM is more appealing to deal with, as it can be efficiently and exactly solved by a convex optimization or conic programming  solver. We further dealt with the joint location  and cost optimization problem under MRUM by proposing two exact methods based on a conic reformulation and the multicut outer-approximation scheme, and a local search heuristic. We provided experiments based on generated instances of various sizes, which show  the efficiency and scalability of our approaches, in particular the CONIC approach for the cost optimization problem and the MOA algorithms for the joint problem. 

Our paper clearly shows several advantages of the MRUM framework in facility location, which would suggest promising directions for future work, for instance, exploring the use of the MRUM framework in other choice-based decision-making problems such as joint assortment and price optimization \citep{wang2012capacitated} or security game problems \citep{yang2012computing}.

\bibliographystyle{plainnat_custom}
\bibliography{refs}

\pagebreak

\appendix
\section*{APPENDIX}

\section{Proof of Proposition \ref{prop:derivative-phi}}
\proof{}
The KKT  conditions imply that, for any $\by\in \cY$,
\begin{align}
    L(\bx^*(\by),\lambda^*(\by),\bgamma^{U*}(\by),\bgamma^{L^*}(\by) |\by) &=  f^{\MRUM}(\by,\bx^{*}(\by)) \label{eq:dual-1}\\
    \frac{L(\bx^*(\by),\lambda^*(\by),\bgamma^{U*}(\by),\bgamma^{L^*}(\by) |\by)}{\partial x_i} &= 0,\; \forall i\in [m]\label{eq:dual-2}\\
    \lambda^*(\by) \left(\sum_{i\in[m]} x^*_i(\by)-C\right) &= 0\label{eq:dual-3}\\
     \gamma^{U*}_i(\by)(x^*_i(\by) - y_i U_i)&= 0,\;\forall i\in [m]\label{eq:dual-4}\\
      \gamma^{L*}_i(\by)(x^*_i(\by) - y_i L_i)&= 0,\;\forall i\in [m]\label{eq:dual-5}
\end{align}
For  notational simplicity, we omit $\by$ from 
$\lambda^*(\by),\bgamma^{U*}(\by)$ and $\bgamma^{L^*}(\by)$ in some places. 
Taking the derivatives of the both sides of \eqref{eq:dual-1} w.r.t. $y_i$, noting that $\Phi(\by) = f^{\MRUM}(\by,\bx^*(\by))$, we have
\begin{align}
    \frac{\Phi(\by)}{\partial y_i}&= \frac{f^{\MRUM}(\by,\bx^*(\by))}{\partial y_i} = \frac{\partial   L(\bx^*(\by),\lambda^*(\by),\bgamma^{U*}(\by),\bgamma^{L*}(\by) |\by)}{\partial y_i } \nonumber\\
    &= \frac{L(\bx^*,\lambda^*,\bgamma^{U*},\bgamma^{L*} |\by)}{\partial y_i } 
    + \sum_{i'\in [m]}\frac{L(\bx^*,\lambda^*,\bgamma^{U*},\bgamma^{L*} |\by)}{\partial x_{i'} }\frac{x^*_{i'}(\by)}{\partial y_i} 
    + \frac{L(\bx^*,\lambda^*,\bgamma^{U*},\bgamma^{L*} |\by)}{\partial \lambda }\frac{\lambda^*(\by)}{\partial y_i} \nonumber \\
    &+ \sum_{i'\in [m]}\frac{L(\bx^*,\lambda^*,\bgamma^{U*},\bgamma^{L*} |\by)}{\partial \gamma_{i'}^{U^*} }\frac{\gamma_{i'}^{U*}(\by)}{\partial y_i} + \sum_{i'\in [m]}\frac{L(\bx^*,\lambda^*,\bgamma^{U*},\bgamma^{L*} |\by)}{\partial \gamma_{i'}^{L^*} }\frac{\gamma_{i'}^{L*}(\by)}{\partial y_i} \label{eq:dual-eq1} 
\end{align}
Now, looking at each term of \eqref{eq:dual-eq1} with the help of the  KKT conditions, we see that 
\begin{equation}\label{eq:cond-1}
 \sum_{i'\in [m]}\frac{L(\bx^*,\lambda^*,\bgamma^{U*},\bgamma^{L*} |\by)}{\partial x_{i'} }\frac{x^*_{i'}(\by)}{\partial y_i}  = 0
\end{equation}
due to \eqref{eq:dual-2}. Moreover, since $\lambda (\sum_{i \in [m]} x_i-C)$ is the only term in $L(\bx,\lambda,\bgamma^U,\bgamma^L |\by)$ involving  $\lambda$, we have 
\[
 \frac{L(\bx^*,\lambda^*,\bgamma^{U*},\bgamma^{L*} |\by)}{\partial \lambda }\frac{\lambda^*(\by)}{\partial y_i}  = \left(\sum_{i'\in [m]} x^*_{j'}(\by)-C\right) \frac{\partial \lambda^*(\by)}{\partial y_i}
\]
If $\sum_{i'\in [m]} x^*_{i'}(\by)-C = 0$, then we have $\frac{L(\bx^*,\lambda^*,\bgamma^{U*},\bgamma^{L*} |\by)}{\partial \lambda }\frac{\lambda^*(\by)}{\partial y_i} = 0$. Otherwise, if $\sum_{i'\in [m]} x^*_{i'}(\by)-C < 0$, the the complementary slackness \eqref{eq:dual-3} tells us that $\lambda^*(\by) = 0$. Moreover, due to the continuity of $\lambda^*(\by)$ and $\bx^*(\by)$ in $\by$, we should have $\frac{\lambda^*(\by)}{\partial y_i} = 0$. Thus, we always have 
\begin{equation}
\label{eq:cond-2}
    \frac{L(\bx^*,\lambda^*,\bgamma^{U*},\bgamma^{L*} |\by)}{\partial \lambda }\frac{\lambda^*(\by)}{\partial y_i} = 0.
\end{equation} 
Similarly,  we also have 
\begin{align}
    \sum_{i'\in [m]}\frac{L(\bx^*,\lambda^*,\bgamma^{U*},\bgamma^{L*} |\by)}{\partial \gamma_{i'}^{U^*} }\frac{\gamma_{i'}^{U*}(\by)}{\partial y_i}  = 0 \label{eq:cond-3} \\
    \sum_{i'\in [m]}\frac{L(\bx^*,\lambda^*,\bgamma^{U*},\bgamma^{L*} |\by)}{\partial \gamma_{i'}^{L^*} }\frac{\gamma_{i'}^{L*}(\by)}{\partial y_i} = 0\label{eq:cond-4}
\end{align}
Combining \eqref{eq:dual-eq1}-\eqref{eq:cond-4}, we have
\begin{align}
\frac{\partial \Phi(\by)}{\partial y_i} &= \frac{\partial   L(\bx^*,\lambda^*,\bgamma^{U*},\bgamma^{L*} |\by)}{\partial y_i } \nonumber\\
&=\frac{\partial f^{\MRUM}(\by,\bx^*) }{\partial y_i} + \gamma^{U*}_i U_i - \gamma^{L*}_i L_i, \;\forall i\in [m],\nonumber
\end{align}
as desired. 
\endproof

\section{Addtional Experiments}

In Table \ref{table3} we report the comparison results for large-sized instances of the joint location and cost optimization problem under MRUM. Each row of the table corresponds to 8 instances, grouped by $(|I|,m, C,K)$, similarly as in Table \ref{table2}. 
There are 36 groups of instances in total.
We also highlight in bold the largest numbers of instances solved to optimality and the largest numbers of instances solved
with the best objectives values. The symbol ``-'' indicates that  the corresponding algorithm encounters an ``out-of-memory'' issue and cannot return any solution. 
In terms of the number of instances solved to optimality, the MOA algorithm with $T=1$ (i.e. M1) seems to perform the best; it gives the largest numbers of instances solved to optimality for 32/36 groups of instances.  Moreover, in these 32/36 groups of instances, M1  solves 8/8 instances  to optimality  for 26 groups. 
Such numbers for the CONIC and  M5 are  8 and 17, respectively, which are remarkably smaller than the M1.
On average,  M1 is able to return optimal solutions for 87.5\% of the instances while the percentages of instances solved to optimality for the  CONIC and  M5 approaches are 50.5\%, 76\%, respectively. 
In terms of number of instances solved with the best objective values, we also observe that while the MOA algorithms  achieve the best objective values for  99.6\% instances, CONIC and LS give best objective values for only  72.8\% and 59\%  of the instances, respectively.

\begin{table}
\centering
\caption{Comparison results for the joint location and cost optimization problem, for\textit{ large-sized} instances grouped by $(m,|I|,C,K)$ (8 instances per row); CN stands for the CONIC approach and M1, M5 stand for the MOA algorithm with $T=1,5$, respectively.}
\label{table3}
\resizebox{\linewidth}{!}{%
\begin{tabular}{lllllll|lll|lll} 
\hline
\multirow{2}{*}{$m$} & \multirow{2}{*}{$|I|$} & \multirow{2}{*}{$C$} & \multirow{2}{*}{$K$} & \multicolumn{3}{c|}{\begin{tabular}[c]{@{}c@{}}\# instances solved\\to optimality\end{tabular}} & \multicolumn{3}{c|}{\begin{tabular}[c]{@{}c@{}}\# instances solved\\with best objectives\end{tabular}} & \multicolumn{3}{c}{Average CPU time (s)} \\ 
\cline{5-13}
 &  &  &  & CN & M1 & M5 & CN & M1/M5 & LS & CN & M1 & M5 \\ 
\hline
1000 & 1000 & 200 & 200 & 1 & \textbf{7 } & \textbf{7 } & \textbf{\textbf{8}} & \textbf{8} & 4 & 539.0 & 208.7 & \textbf{189.4} \\ 
\hline
1000 & 1000 & 200 & 500 & 4 & \textbf{8 } & \textbf{8 } & \textbf{\textbf{8}} & \textbf{8} & 0 & 382.0 & \textbf{8.0} & 13.7 \\ 
\hline
1000 & 1000 & 200 & 700 & 7 & \textbf{8 } & \textbf{8 } & \textbf{\textbf{8}} & \textbf{8} & 0 & 271.5 & 3.5 & \textbf{3.0} \\ 
\hline
1000 & 1000 & 500 & 200 & 0 & \textbf{6 } & \textbf{6 } & \textbf{\textbf{8}} & \textbf{8} & 4 & 600.0 & 329.5 & \textbf{227.3} \\ 
\hline
1000 & 1000 & 500 & 500 & 7 & \textbf{8 } & \textbf{8 } & \textbf{\textbf{8}} & \textbf{8} & 4 & 281.2 & \textbf{15.9} & 62.1 \\ 
\hline
1000 & 1000 & 500 & 700 & 6 & \textbf{8 } & \textbf{8 } & \textbf{\textbf{8}} & \textbf{8} & 4 & 338.6 & \textbf{6.8} & 7.9 \\ 
\hline
1000 & 1000 & 700 & 200 & 0 & \textbf{6 } & 4 & \textbf{\textbf{8}} & \textbf{8} & \textbf{\textbf{8}} & 600.0 & 404.7 & \textbf{221.1} \\ 
\hline
1000 & 1000 & 700 & 500 & 7 & \textbf{8 } & \textbf{8 } & \textbf{\textbf{8}} & \textbf{8} & 4 & 282.2 & \textbf{38.7} & 87.9 \\ 
\hline
1000 & 1000 & 700 & 700 & 7 & \textbf{8 } & \textbf{8 } & \textbf{\textbf{8}} & \textbf{8} & 4 & 229.5 & \textbf{94.0} & 97.4 \\ 
\hline
1000 & 5000 & 200 & 200 & 0 & 0 & 0 & 0 & \textbf{8} & 0 & \textbf{600.0} & \textbf{600.0} & \textbf{600.0} \\ 
\hline
1000 & 5000 & 200 & 500 & 0 & \textbf{8 } & \textbf{8 } & 0 & \textbf{8} & 0 & 600.0 & \textbf{87.0} & 98.0 \\ 
\hline
1000 & 5000 & 200 & 700 & 0 & \textbf{8 } & \textbf{8 } & 0 & \textbf{8} & 0 & 600.0 & 9.7 & \textbf{4.8} \\ 
\hline
1000 & 5000 & 500 & 200 & 0 & 0 & 0 & 0 & \textbf{8} & 4 & \textbf{600.0} & \textbf{600.0} & \textbf{600.0} \\ 
\hline
1000 & 5000 & 500 & 500 & 0 & \textbf{8 } & 7 & 0 & \textbf{8} & 4 & 600.0 & \textbf{268.2} & 328.8 \\ 
\hline
1000 & 5000 & 500 & 700 & 0 & \textbf{7 } & \textbf{7 } & 0 & \textbf{8} & 4 & 600.0 & 164.4 & \textbf{116.7} \\ 
\hline
1000 & 5000 & 700 & 200 & 0 & 0 & 0 & 0 & \textbf{8} & \textbf{\textbf{8}} & \textbf{600.0} & \textbf{600.0} & \textbf{600.0} \\ 
\hline
1000 & 5000 & 700 & 500 & 0 & \textbf{5 } & 4 & 0 & \textbf{8} & 4 & 600.0 & 429.5 & \textbf{400.0} \\ 
\hline
1000 & 5000 & 700 & 700 & 0 & 6 & \textbf{7 } & 0 & \textbf{8} & 4 & 600.0 & 158.0 & \textbf{7.7} \\ 
\hline
2000 & 200 & 400 & 400 & 4 & \textbf{8 } & 7 & \textbf{\textbf{8}} & \textbf{8} & 4 & 378.5 & 43.7 & \textbf{39.3} \\ 
\hline
2000 & 200 & 400 & 1000 & \textbf{8 } & \textbf{8 } & \textbf{8 } & \textbf{\textbf{8}} & \textbf{8} & 4 & 21.0 & \textbf{5.0} & 10.4 \\ 
\hline
2000 & 200 & 400 & 1400 & \textbf{8 } & \textbf{8 } & \textbf{8 } & \textbf{\textbf{8}} & \textbf{8} & 4 & 16.9 & \textbf{5.8} & 7.2 \\ 
\hline
2000 & 200 & 1000 & 400 & 6 & \textbf{8 } & 5 & \textbf{\textbf{8}} & \textbf{8} & \textbf{\textbf{8}} & 240.1 & 112.8 & \textbf{83.3} \\ 
\hline
2000 & 200 & 1000 & 1000 & 6 & \textbf{8 } & \textbf{8 } & \textbf{\textbf{8}} & \textbf{8} & \textbf{\textbf{8}} & 169.6 & \textbf{8.7} & 21.1 \\ 
\hline
2000 & 200 & 1000 & 1400 & \textbf{8 } & \textbf{8 } & 7 & \textbf{\textbf{8}} & \textbf{8} & 4 & 23.4 & \textbf{13.9} & 16.1 \\ 
\hline
2000 & 200 & 1400 & 400 & 6 & \textbf{8 } & 4 & \textbf{\textbf{8}} & \textbf{8} & \textbf{\textbf{8}} & 189.5 & 107.7 & \textbf{55.1} \\ 
\hline
2000 & 200 & 1400 & 1000 & \textbf{8 } & \textbf{8 } & 7 & \textbf{\textbf{8}} & \textbf{8} & \textbf{\textbf{8}} & 30.8 & \textbf{11.2} & 23.5 \\ 
\hline
2000 & 200 & 1400 & 1400 & \textbf{8 } & \textbf{8 } & 7 & \textbf{\textbf{8}} & \textbf{8} & \textbf{\textbf{8}} & 33.0 & 14.7 & \textbf{12.2} \\ 
\hline
3000 & 100 & 600 & 600 & 5 & \textbf{8 } & 5 & \textbf{\textbf{8}} & \textbf{8} & 4 & 326.2 & \textbf{23.8} & 54.1 \\ 
\hline
3000 & 100 & 600 & 1500 & 7 & \textbf{8 } & \textbf{8 } & \textbf{\textbf{8}} & \textbf{8} & 4 & 90.0 & \textbf{4.9} & 15.6 \\ 
\hline
3000 & 100 & 600 & 2100 & \textbf{8 } & \textbf{8 } & 7 & \textbf{\textbf{8}} & \textbf{8} & 4 & 24.7 & \textbf{7.8} & 11.2 \\ 
\hline
3000 & 100 & 1500 & 600 & 6 & \textbf{7 } & 2 & \textbf{\textbf{8}} & 7 & 7 & 173.9 & \textbf{31.5} & 49.6 \\ 
\hline
3000 & 100 & 1500 & 1500 & 5 & \textbf{8 } & \textbf{8 } & 6 & \textbf{8} & 7 & 237.3 & \textbf{9.0} & 25.2 \\ 
\hline
3000 & 100 & 1500 & 2100 & \textbf{8 } & \textbf{8 } & 7 & \textbf{\textbf{8}} & \textbf{8} & 7 & 29.3 & \textbf{10.9} & 15.6 \\ 
\hline
3000 & 100 & 2100 & 600 & 5 & \textbf{8 } & 2 & \textbf{\textbf{8}} & \textbf{8} & 7 & 239.5 & 49.2 & \textbf{14.1} \\ 
\hline
3000 & 100 & 2100 & 1500 & 4 & \textbf{8 } & 6 & 4 & \textbf{8} & 7 & 312.3 & 15.1 & \textbf{13.2} \\ 
\hline
3000 & 100 & 2100 & 2100 & \textbf{8 } & \textbf{8 } & 7 & \textbf{\textbf{8}} & \textbf{8} & 7 & 66.5 & \textbf{7.5} & 14.7 \\ 
\hline
\multicolumn{4}{l}{Average} & 4.36 & 7 & 6.08 & 5.83 & 7.97 & 4.72 &  &  &  \\
\hline
\end{tabular}
}
\end{table}

In terms of average CPU time, we first note that LS remains to be the worst since it exceeds the time budget of 600 seconds for all the instances. In fact, for  medium-size and large-size instances, due to the expensiveness of the computation of $\Phi(\by)$, 
the LS algorithm could not complete the first step (i.e., adding locations to an empty set)  within the time budget.
Moreover,  if we look at the groups of instances that all the three approaches (CONIC, M1, M5) provide  optimal solutions for all 8/8 instances, the average CPU times required by CONIC are higher than those required by M1 and M5.
We also observe that CONIC  performs well for instances of  $|I|\leq 200$, $m \geq 1000$ and $K = 0.7m$. For these instances, CONIC is even better than M5. In summary, MOA with $T =1$ requires the shortest CPU times in most groups of instances (22/36). It is better than the MOA algorithm with $T =5$ (the numbers of such groups for M5 is 16), although the gaps are not as large as in the small-size and medium-size instances. The CONIC and LS approaches always require higher average CPU times for all the groups of instances,  as compared to the MOA algorithms.

\end{document}